\documentclass[a4,11pt]{amsart}
\usepackage{amssymb,amsmath,amsthm,txfonts}
\usepackage{cite}

\usepackage{color}

\newtheorem{thm}{Theorem}[section]
\newtheorem{lem}[thm]{Lemma}

\newtheorem{prop}[thm]{Proposition}
\theoremstyle{definition}

\theoremstyle{remark}
\newtheorem{rem}[thm]{\textbf{Remark}}
\newtheorem{rems}[thm]{\textbf{Remarks}}

 \makeatletter
    
    \@addtoreset{equation}{section}
  \makeatother

\makeatletter
      \def\@makefnmark{%
         \leavevmode
            \raise.9ex\hbox{\check@mathfonts
                \fontsize\sf@size\z@\normalfont%
                            \@thefnmark}%
       }
      \makeatother

\newcommand{\D}{\textrm{div}}

\newcommand{\dd}{\textrm{d}}

\begin{document}

\title[]{Stability of Lamb dipoles}
\author[]{K.Abe and K.Choi}
\date{}
\address[K. Abe]{Department of Mathematics, Graduate School of Science, Osaka City University, 3-3-138 Sugimoto, Sumiyoshi-ku Osaka, 558-8585, Japan}
\email{kabe@sci.osaka-cu.ac.jp}
\address[K. Choi]{Department of Mathematical Sciences, Ulsan National Institute of Science and Technology, UNIST-gil 50, Ulsan, 44919, Republic of Korea}
\email{kchoi@unist.ac.kr}

\subjclass[2010]{35Q35, 35K90}
\keywords{Lamb dipole, Euler equations, orbital stability, vortex pairs}
\date{\today}

\maketitle


\begin{abstract}
The Lamb dipole is a traveling wave solution to the two-dimensional Euler equations introduced by S. A. Chaplygin (1903) and H. Lamb (1906) at the early 20th century. We prove orbital stability of this solution based on a vorticity method initiated by V. I. Arnold. Our method is a minimization of a penalized energy with multiple constraints that deduces existence and orbital stability for a family of traveling waves. As a typical case, orbital stability of the Lamb dipole is deduced by characterizing a set of minimizers as an orbit of the dipole by a uniqueness theorem in the variational setting.
\end{abstract}

\vspace{15pt}

\section{Introduction}

\subsection{Lamb dipoles}

\vspace{15pt}

We consider the two-dimensional vorticity equations:

\begin{equation*}
\begin{aligned}
\partial_t \zeta+v\cdot \nabla \zeta=0,\quad  v&=k*\zeta \quad \textrm{in}\ \mathbb{R}^{2}\times (0,\infty),\\
\zeta&=\zeta_0\hspace{22pt} \textrm{on}\ \mathbb{R}^{2}\times \{t=0\},
\end{aligned}
\tag{1.1}
\end{equation*}\\ 
with the kernel $k(x)=(2\pi)^{-1}x^{\perp}|x|^{-2}$, $x^{\perp}={}^{t}(-x_2,x_1)$. The equations (1.1) admit \textit{a vortex pair}, i.e., a solution of the form

\begin{equation*}
\begin{aligned}
v(x,t)&=u(x+u_{\infty}t)-u_{\infty},\\
\zeta(x,t)&=\omega(x+u_{\infty}t),    
\end{aligned}
\end{equation*}\\
vanishing at space infinity with a constant velocity $u_\infty\in \mathbb{R}^{2}$. Vortex pairs are symmetric dipoles with compactly supported two vorticities having opposite signs translating to one direction. They are theoretical models of coherent vortex structures in large-scale geophysical flows. See, e.g., \cite{VF89}, \cite{FV94} for experimental works. By rotational invariance of (1.1), we take $u_{\infty}={}^{t}(-W,0)$, $W>0$, without loss of generality. Substituting $(v,\zeta)$ into (1.1) implies the steady Euler equations for $(u,\omega)$ in a half plane:

\begin{equation*}
\begin{aligned}
u\cdot \nabla \omega&=0\qquad \textrm{in}\ \mathbb{R}^{2}_{+},\\ 
u&\to u_{\infty}\quad \textrm{as}\ |x|\to\infty.
\end{aligned}
\tag{1.2}
\end{equation*}

In 1906, H. Lamb \cite[p.231]{Lamb} noted an explicit solution to (1.2), generally referred to as \textit{the Lamb dipole} (Chaplygin-Lamb dipole), a solution $\omega_{L}=\lambda \max\{\Psi_{L},0 \}$, $u_{L}={}^{t}(\partial_{x_2}\Psi_{L},-\partial_{x_1}\Psi_{L})$, $0<\lambda<\infty$, of the form

\begin{equation*}
\Psi_{L}(x)=\left\{
\begin{aligned}
C_LJ_1(\lambda^{1/2}r)\sin\theta,\quad r\leq a,\\
-W\left(r-\frac{a^{2}}{r}\right)\sin\theta,\quad r>a,
\end{aligned}
\right. \tag{1.3}
\end{equation*}\\
with the constants

\begin{align*}
C_L=-\frac{2W}{\lambda^{1/2}J_0(c_0)},\quad a=c_0\lambda^{-1/2},   
\end{align*}\\
where $(r,\theta)$ is the polar coordinate and $J_{m}(r)$ is the $m$-th order Bessel function of the first kind. The constant $c_0$ is the first zero point of $J_1$, i.e., $J_1(c_0)=0$, $c_0=3.8317\cdots$, $J_0(c_0)<0$. The parameter $\lambda>0$ denotes the strength of the vortex and is related with its impulse by

\begin{align*}
\int_{\mathbb{R}^{2}_{+}}x_2\omega_{L}\dd x=\frac{c_0^{2}\pi W}{\lambda}.
\end{align*}\\
The Lamb dipole (1.3) is the simplest explicit solution to (1.2), symmetric for the $x_2$-variable, which is a special case of non-symmetric Chaplygin dipoles, independently founded by S. A. Chaplygin in 1903 \cite{Chap1903}, \cite{Chap07}. See also \cite{MV94}. 

The Lamb dipole is considered as a stable vortex structure in a two-dimensional flow. Its stability has been studied by an experimental work \cite{FV94} and also by a numerical work \cite{VV98}. On the other hand, despite the explicit form of this classical solution, its mathematical stability had been an open question since the solution was introduced by S. A. Chaplygin and H. Lamb at the early 20th century. For solutions with a single-signed vortex such as a circular vortex \cite{WP85}, \cite{SV09} or a rectangular vortex \cite{Denisov17}, stability results have been developed, while no stability result was known for the Lamb dipole which has a multi-signed vortex and forms a traveling wave.

There is an interesting relation with \textit{solitons} in the theory of nonlinear wave equations. One of classical models that describes propagation of a wave may be the KdV equation \cite{KdV}. More generally for the gKdV equation, 

\begin{align*}
\partial_t w+\partial_{x}^{3}w+\partial_{x}(w^{p})=0,\quad x\in \mathbb{R},\ t>0,
\end{align*}\\
for an integer $p\geq 2$, there exists a soliton solution of the form $w(x,t)=Q_{c}(x-ct)$ for $c>0$ and $Q_{c}(x)=c^{1/(p-1)} Q(c^{1/2}x)$, where 

\begin{align*}
Q(x)=\left( \frac{p+1}{2\cosh^{2}((p-1)x /2)}     \right)^{1/(p-1)},
\end{align*}\\
is called soliton, which is a unique positive solution of the elliptic problem $\partial_x^{2}Q+Q^{p}=Q$, up to translation. Stability of this soliton is well known when the problem is globally well-posed. Indeed for $2\leq p<5$, the gKdV equation is globally well-posed, and if initial data is close to the soliton, the solution remains nearby the soliton for all time by admitting translation of $Q$ \cite{Ben72}, \cite{We86}. Such stability is termed \textit{orbital stability}. For $p=5$, this soliton is unstable \cite{MartelMerle01} and a finite time blow-up occurs \cite{Merle01}, \cite{MartelMerle02}. The Euler equations may have some aspects of the wave equation. Even for the three-dimensional case, vortex rings form traveling waves. We shall establish the orbital stability theorem for the Lamb dipole which is the most typical traveling wave.

In the sequel, we identify a function $\zeta_0$ in $\mathbb{R}^{2}_{+}$ with an odd extension to $\mathbb{R}^{2}$ for the $x_2$-variable, i.e., $\zeta_0(x_1,x_2)=-\zeta_0(x_1,-x_2)$. Since a classical solution to (1.1) exists and is symmetric for the $x_2$-variable for sufficiently smooth initial data \cite{LNX01}, a standard approximation argument implies the existence of a symmetric global weak solution $\zeta\in BC([0,\infty); L^{2}\cap L^{1}(\mathbb{R}^{2}) )$ for symmetric initial data $\zeta_{0}\in L^{2}\cap L^{1}(\mathbb{R}^{2})$ \cite{MaB}. Here, $BC([0,\infty); X)$ denotes the space of all bounded continuous functions from $[0,\infty)$ into a Banach space $X$. Among other results, our simplest result is the following:

\vspace{15pt}

\begin{thm}
Let $0<\lambda, W<\infty$. The Lamb dipole $\omega_{L}$ is orbitally stable in the sense that for $\nu>0$ and $\varepsilon>0$, there exists $\delta>0$ such that for $\zeta_0\in L^{2}\cap L^{1}(\mathbb{R}^{2}_{+})$ satisfying $x_2\zeta_{0}\in L^{1}(\mathbb{R}^{2}_{+})$, $\zeta_{0}\geq 0$, $||\zeta_0||_1\leq \nu$ and 

\begin{align*}
\inf_{y\in \partial\mathbb{R}^{2}_{+}}\left\{\left\|\zeta_0-\omega_{L}(\cdot+y) \right\|_{2}+\left\|x_2(\zeta_0-\omega_{L}(\cdot+y)) \right\|_{1}\right\}\leq \delta,
\end{align*}\\
there exists a global weak solution $\zeta(t)$ of (1.1) satisfying

\begin{align*}
\inf_{y\in \partial\mathbb{R}^{2}_{+}}\left\{\left\|\zeta(t)-\omega_{L}(\cdot+y) \right\|_{2}+\left\|x_2(\zeta(t)-\omega_{L}(\cdot+y)) \right\|_{1}\right\}\leq \varepsilon,\quad \textrm{for all}\ t\geq0.  
\end{align*} 
\end{thm}

\vspace{15pt}

\begin{rem}
As we will see later in Remarks 5.2 (i), the smallness condition in Theorem 1.1 can be replaced with a slightly weaker condition $\inf_{y\in \partial\mathbb{R}^{2}_{+}}\left\|\zeta_0-\omega_{L}(\cdot+y) \right\|_{2}+\left|\int x_2\zeta_0\dd x-\mu \right|$ $\leq \delta$ for $\mu=c_0^{2}\pi W/\lambda$. 
\end{rem}

\vspace{15pt}

\subsection{Vorticity method}

Theorem 1.1 is a particular case of our general stability theorem. Let us consider the existence problem (1.2). The equation $(1.2)_1$ can be written by the Jacobian of ${}^{t}(\Psi,\omega)$ for $u={}^{t}(\partial_{x_2}\Psi,-\partial_{x_1}\Psi )$. Therefore $\omega$ is represented by $\omega=\lambda f(\Psi)$ with some function $f(t)$ and $\lambda>0$ and existence of such $(u,\omega)$ is reduced to the free-boundary problem for $\gamma\geq 0$:

\begin{equation*}
\left.
\begin{aligned}
-\Delta\Psi&=\lambda f(\Psi)\quad \textrm{in}\ \mathbb{R}^{2}_{+},\\
\Psi&=-\gamma\hspace{25pt} \textrm{on}\ \partial\mathbb{R}^{2}_{+},\\
\partial_{x_1}\Psi&\to0,\quad \partial_{x_2}\Psi\to-W\quad \textrm{as}\ |x|\to \infty.
\end{aligned}
\right.\tag{1.4}
\end{equation*}\\
The function $f$ is called a vorticity function which is prescribed by a non-negative and non-decreasing function. In this paper, we shall take

\begin{align*}
f(t)=t_{+},\quad t_{+}=\max\{t,0\},
\end{align*}\\
for which the Lamb dipole $\Psi_{L}$ is a solution to (1.4) for $\gamma=0$ and $\textrm{spt}\ \omega_{L}=\overline{B(0,a)\cap \mathbb{R}^{2}_{+}}$, i.e., $\omega_{L}=\lambda f(\Psi_{L})$. Here $B(0,a)$ is an open disk centered at the origin with the radius $a>0$. The three parameters $W, \gamma\geq 0$ and $\lambda>0$ are referred to as propagation speed, flux constant and strength parameter. We chose the flux constant $\gamma$ so that $\Psi=0$ on the boundary of the vortex core $\textrm{spt}\ \omega=\overline{\Omega}$. The problem (1.4) is a free-boundary problem since the vortex core $\Omega$ is a priori unknown. Once the core is found, one can find $\Psi$ by solving the two problems: 

\begin{align*}
&-\Delta \Psi=\lambda \Psi \quad \textrm{in}\ \Omega,\quad \Psi=0\quad \textrm{on}\ \partial\Omega, \\
&-\Delta\Psi=0\quad \textrm{in}\ \mathbb{R}^{2}_{+}\backslash \overline{\Omega},\ \Psi=-\gamma\quad \textrm{on}\ \partial\mathbb{R}^{2}_{+},\ \partial_{x_1}\Psi\to0,\ \partial_{x_2}\Psi\to -W\quad \textrm{as}\ |x|\to\infty.
\end{align*}\\
On the other hand, the core is characterized as $\Omega=\{x\in \mathbb{R}^{2}_{+}\ |\ \Psi(x)>0 \}$ by a maximum principle. The function $\Psi=\psi-Wx_2-\gamma$ is represented by the Green function of the Laplace operator subject to the Dirichlet boundary condition in a half plane

\begin{align*}
\psi(x)=\int_{\mathbb{R}^{2}_{+}}G(x,y)\omega(y)\dd y, \quad G(x,y)=\frac{1}{4\pi}\log{\left(1+\frac{4x_2y_2}{|x-y|^{2}} \right)}.  \tag{1.5}
\end{align*}

To study existence and stability of solutions to (1.4), we consider a variational principle based on vorticity, called a vorticity method, originating from the idea of Kelvin \cite{kelvin1880}, initiated by Arnold \cite{Arnold66}, \cite{AK98}. See also Benjamin \cite{Ben76} for vortex rings. For vortex pairs, vorticity methods were developed by Turkington \cite{Tu83} and Burton \cite{Burton88}. See also Norbury \cite{Norbury75} and Yang \cite{Yang91} for a stream function method.

Our approach is based on the vorticity method of Friedman-Turkington \cite{FT81}, \cite{Friedman82} developed for vortex rings. For $0<\mu,\nu,\lambda<\infty$, we set a space of admissible functions

\begin{align*}
K_{\mu,\nu}=\left\{\omega\in L^{2}(\mathbb{R}^{2}_{+})\ \middle|\ \omega\geq 0,\ \int_{\mathbb{R}^{2}_{+}} x_2\omega\dd x=\mu,\ \int_{\mathbb{R}^{2}_{+}}\omega\dd x\leq \nu\    \right\}.
\end{align*}\\
We construct solutions of (1.4) by maximizing a penalized energy 

\begin{align*}
E_{2,\lambda}[\omega]=E[\omega]-\frac{1}{2\lambda}\int_{\mathbb{R}^{2}_{+} }\omega^{2}\dd x, \quad E[\omega]=\frac{1}{2}\int_{\mathbb{R}^{2}_{+}}\int_{\mathbb{R}^{2}_{+}}G(x,y)\omega(x)\omega(y)\dd x\dd y.
\end{align*}\\
For a notational convenience, we formulate the maximization problem as a minimization of $-E_{2,\lambda}$ and denote by

\begin{align*}
I_{\mu,\nu\,\lambda}=\inf_{\omega\in K_{\mu,\nu}}\left\{-E_{2,\lambda}[\omega]\right\}.   \tag{1.6}
\end{align*}\\
The constants $W,\ \gamma\geq 0$ are Lagrange multipliers. This formulation is slightly different from that of \cite{FT81}, \cite{Friedman82}, where admissible functions are restricted to a space of symmetric functions for $x_1\in \mathbb{R}$. More precisely, the method in \cite{FT81}, \cite{Friedman82} applies to prove compactness of a minimizing sequence satisfying 

\begin{equation*}
\begin{aligned}
&\omega(x_1,x_2)=\omega(-x_1,x_2), \\
&\omega(x_1,x_2)\ \textrm{is non-increasing for}\ x_1>0.
\end{aligned}
\tag{1.7}
\end{equation*}\\
The condition (1.7) is essential for the method in \cite{FT81}, \cite{Friedman82}. In fact, since the energy $-E_{2,\lambda}$ is invariant by translation for the $x_1$-variable, translation of any minimizer is a minimizing sequence. In this paper, without assuming (1.7), we shall show that any minimizing sequence is relatively compact by translation for the $x_1$-variable by using the concentration compactness principle of Lions \cite{Lions84a}. The following Theorem 1.3 is an improvement of \cite{FT81}, \cite{Friedman82} in terms of vortex pairs. 
 
\vspace{15pt}

\begin{thm}
Let $0<\mu,\nu, \lambda<\infty$. For any minimizing sequence $\{\omega_n\}$ satisfying $\omega_n \in K_{\mu_n,\nu}$, $\mu_n\to \mu$ and $-E_{2,\lambda}[\omega_n]\to I_{\mu,\nu,\lambda}$, there exists a sequence $\{y_{n}\}\subset \partial\mathbb{R}^{2}_{+}$ such that $\{\omega_n(\cdot+y_n)\}$ and $\{x_2\omega_n(\cdot+y_n)\}$ are relatively compact in $L^{2}(\mathbb{R}^{2}_{+})$ and $L^{1}(\mathbb{R}^{2}_{+})$, respectively. In particular, the problem (1.6) has a minimizer in $K_{\mu,\nu}$.  
\end{thm}

\vspace{15pt}
There is a novelty to adapt the vorticity method of \cite{FT81}, \cite{Friedman82}, instead of \cite{Tu83} which prescribes that mass is exactly $\nu>0$ for admissible functions. As proved in \cite{FT81}, \cite{Friedman82} for vortex rings, mass becomes strictly less than $\nu>0$ for small $\lambda>0$ with fixed $\mu,\nu$. Indeed, the variational principle in \cite{Tu83} does not provide solutions of (1.4) for small $\lambda>0$. Our existence for small $\lambda>0$ seems a new result although the above formulation is noted in \cite{Tu83}. See also \cite{Norbury75}.

Removing the restriction on the strength parameter is essential in the present work since solutions of (1.4) approach a Lamb dipole as $\lambda\to0$. We shall rigorously state this claim as in Theorem 1.5 below. For fixed $\mu, \nu$, solutions of (1.6) form one parameter family for $0<\lambda<\infty$. In particular, solutions approach a Dirac measure as $\lambda\to\infty$ and in contrast a Lamb dipole as $\lambda\to0$. A variational characterization of the Lamb dipole is studied in \cite{Burton96}, \cite{Burton05b} for solutions to (1.4) for $\gamma=0$.

Orbital stability of vortex pairs is a consequence of compactness of a minimizing sequence. We use conservations of $L^{q}$-norms, impulse and penalized energy of (1.1):

\begin{equation*}
\begin{aligned}
||\zeta||_{q}(t)&=||\zeta_0||_{q},\quad 1\leq q\leq 2,  \\
||x_2\zeta||_{1}(t)&=||x_2\zeta_0||_{1},  \\
E_{2,\lambda}(\zeta)(t)&=E_{2,\lambda}(\zeta_0),\qquad \textrm{for all}\ t\geq 0.
\end{aligned}
\tag{1.8}
\end{equation*}\\
Although a global weak solution $\zeta(t)$ of (1.1) obtained by an approximation argument \cite{MaB} might have weak regularity at $t=0$, by the renormalization property of DiPerna-Lions \cite{DL89}, the constructed weak solution satisfies the conservations (1.8), i.e., $\zeta(t)\in K_{\mu,\nu}$ for $\zeta_0\in K_{\mu,\nu}$. In general, $\zeta(t)$ is not symmetric and non-increasing for the $x_1$-variable even if $\zeta_0$ is. 

The vorticity method not only constructs stationary solutions as lowest energy solutions but also deduces their stability by compactness of a minimizing sequence, cf. \cite{CL82} for dispersive equations. For the Euler equations, research on orbital stability goes back to Benjamin \cite{Ben76}. See Wan \cite{Wan88} for an early work. For vortex pairs, the first orbital stability result appeared in Burton, Nussenzveig Lopes and Lopes Filho \cite{BNL13} for a certain class of solutions to (1.2) by a vorticity method based on a rearrangement for a prescribed function. See \cite{ILN03}, \cite{BNL13} for a physical background and an introduction to the problem. The method of \cite{BNL13} yields existence of solutions to (1.4) for small $W>0$, $\gamma=0$ with unknown $f(t)$, $\lambda>0$ and deduces their stability for compactly supported $\zeta_0$. We prove existence of (1.4) by prescribing $f(t)=t_{+}$, $\lambda>0$ and deduce their stability without assuming compact support for $\zeta_0$. Let $S_{\mu,\nu,\lambda}$ denote the set of minimizers of (1.6). Theorem 1.3 implies: 

\vspace{15pt}

\begin{thm}
For $0<\mu,\nu,\lambda<\infty$, $S_{\mu,\nu,\lambda}$ is orbitally stable in the sense that for $\varepsilon>0$, there exists $\delta>0$ such that for $\zeta_0\in L^{2}\cap L^{1}(\mathbb{R}^{2}_{+})$ satisfying $x_2\zeta_{0}\in L^{1}(\mathbb{R}^{2}_{+})$, $\zeta_0\geq 0$, $||\zeta_0||_1\leq \nu$ and 

\begin{align*}
\inf_{\omega\in S_{\mu,\nu,\lambda}}\left\{||\zeta_0-\omega||_{2}
+||x_2(\zeta_0-\omega)||_{1}\right\}
\leq \delta,
 \tag{1.9}
\end{align*}\\
there exists a global weak solution $\zeta(t)$ of (1.1) satisfying 

\begin{align*}
\inf_{\omega\in S_{\mu,\nu,\lambda}}\left\{||\zeta(t)-\omega||_{2}
+||x_2(\zeta(t)-\omega)||_{1}\right\}
\leq \varepsilon,\quad \textrm{for all}\ t\geq0.   \tag{1.10}
\end{align*}
\end{thm}

\vspace{15pt}
Theorem 1.4 is a general stability theorem for a family of vortex pairs for $0<\lambda<\infty$. If the set of minimizers is characterized as an orbit ${\mathcal{O}}(\omega)=\{\omega(\cdot+y)\ |\ y\in \partial \mathbb{R}^{2}_{+} \}$ for some vortex pair, one can deduce orbital stability of the vortex pair itself. Since translation of a minimizer $\omega$ of (1.6) is also a minimizer, the orbit ${\mathcal{O}}(\omega)$ is a subset of $S_{\mu,\nu,\lambda}$. The converse inclusion is a uniqueness issue. See \cite{AF86} for uniqueness of the Hill's spherical vortex rings and \cite{Burton96}, \cite{Burton05b} of the Lamb dipoles.

In this paper, we prove uniqueness of minimizers of (1.6) for small $\lambda>0$, i.e., $\mu\nu^{-1}\lambda^{1/2}\leq M_1$ for some $M_1>0$. As proved later, the flux constant $\gamma$ vanishes for small $\lambda>0$ and $\psi/x_2$ is a positive solution of the elliptic problem in $\mathbb{R}^{4}$, i.e., for $y={}^{t}(y',y_4)\in \mathbb{R}^{4}$,

\begin{align*}
-\Delta_y \left(\frac{\psi(y_4,|y'|)}{|y'|}  \right)=\lambda f\left(\frac{\psi(y_4,|y'|)}{|y'|}-W\right)\qquad \textrm{in}\ \mathbb{R}^{4}.
\end{align*}\\
Since positive solutions $\psi/|y'|$ of the above problem are radially symmetric for some point on $\{y'=0\}$ \cite{Burton96}, minimizers of (1.6) for small $\lambda>0$ must be translation of a Lamb dipole $\omega_{L}$ for $W>0$. As a consequence, it turns out that $S_{\mu,\nu,\lambda}={\mathcal{O}}(\omega_{L})$ for $\mu\nu^{-1}\lambda^{1/2}\leq M_1$ and (1.10) is orbital stability of the Lamb dipole itself. By the constraint on the impulse, the speed $W>0$ is uniquely determined by $W=\mu\lambda/(c_0^{2}\pi)$.

\vspace{15pt}

\begin{thm}
Let $0<\mu,\nu, \lambda<\infty$ satisfy $\mu\nu^{-1}\lambda^{1/2} \leq M_1$ for some absolute constant $M_1>0$. Let $\omega_L$ be the Lamb dipole for $W=\mu\lambda/(c_0^{2}\pi)$. Then, minimizers of (1.6) are  translation of the Lamb dipole, i.e.,

\begin{align*}
S_{\mu,\nu,\lambda}=\left\{\omega_{L}(\cdot+y)\ \middle|\ y\in \partial\mathbb{R}^{2}_{+} \right\}.    \tag{1.11}
\end{align*}
\end{thm}

\vspace{15pt}

The characterization (1.11) implies that $S_{\mu,\nu,\lambda}$ is independent of large $\nu>0$ for fixed $\mu,\lambda$, i.e., $\mu\nu^{-1}\lambda^{1/2}\leq M_1$. Therefore for given $\lambda,W>0$, $\nu>0$ and $\mu=c_0^{2}\pi W/\lambda$, we take $\tilde{\nu}=\max\{\nu, \mu \lambda^{1/2}M_1^{-1}\}$ so that $S_{\mu,\tilde{\nu},\lambda}={\mathcal{O}}(\omega_L)$. Theorem 1.1 is then deduced from Theorem 1.4. 

There is a possibility that uniqueness still holds for solutions to (1.4) for small $\gamma>0$. See \cite{Norbury72}, \cite{AF88} for uniqueness of vortex rings. If the uniqueness holds, one can characterize $S_{\mu,\nu,\lambda}$ as an orbit of some deformed vortex pair supported away from the boundary $\partial\mathbb{R}^{2}_{+}$. Theorem 1.4 may include stability of such solutions. 

There are few remarks related with nonlinear wave equations. Orbital stability is concerned with stability about a shape of a wave. Indeed, Theorem 1.1 implies that the shape of $\omega_{L}$ is stable by a perturbation for all $t\geq 0$. A more advanced question is the asymptotic behavior of the perturbation $\zeta(t)$ as $t\to\infty$. One may expect that a perturbation approaches some fixed traveling wave as $t\to\infty$. Such stability is termed \textit{asymptotic stability} in the study of nonlinear wave equations. Another issue is interaction between traveling waves. Stability of two Lamb dipoles or more generally stability of a finite number of the dipoles are questions. We refer to a survey \cite{Tao09} on stability of solitons.

In this paper, we considered the vorticity function $f(t)=t_{+}$ to prove the orbital stability of the Lamb dipole. Our method is also applied to prove orbital stability of more general vortex pairs and also vortex rings. For example, we are able to take $f(t)=t^{1/(p-1)}_{+}$ as a vorticity function to study existence and orbital stability of vortex pairs for $4/3<p<\infty$ and vortex rings for $6/5<p<\infty$. The stability norm can be replaced with the $L^{p}$-norm with the weighted $L^{1}$-norm. 

A special case is $p=\infty$ for which the vorticity function becomes an indicator function. The penalized energy can be replaced with the kinetic energy whose minimizers are vortex patches \cite{FT81}, \cite{Friedman82}. In contrast to the stability of the circular vortex \cite{WP85}, \cite{SV09}, orbital stability of translating patches are questions. This class particularly includes the Hill's spherical vortex rings.

\vspace{15pt}

Let us sketch the proof of Theorem 1.3. In the sequel, we reduce the problem to the case $\nu=\lambda=1$ by the scaling 

\begin{align*}
\hat{\omega}(x)=\frac{1}{\lambda \nu}\omega\left(\frac{x}{\lambda^{1/2}}\right).    \tag{1.12}
\end{align*}\\
If $\omega\in K_{\mu,\nu}$, $\hat{\omega}\in K_{M,1}$ for $M=\mu\nu^{-1}\lambda^{1/2}$ and $E_{2,1}[\hat{\omega}]=\nu^{-2}E_{2,\lambda}[\omega]$. We abbreviate the notation as $K_{\mu}=K_{\mu,1}$, $I_{\mu}=I_{\mu,1,1}$, $E_{2}[\omega]=E_{2,1}[\omega]$, and $S_{\mu}=S_{\mu,1,1}$. 

To prove compactness of a minimizing sequence of (1.6), we apply a concentration compactness principle and exclude possibilities of dichotomy and vanishing of the sequence. Since $I_{\mu}$ is negative and decreasing for $\mu\in (0,\infty)$, vanishing can not occur. The problem is to exclude dichotomy of the sequence. Let us consider for simplicity a minimizing sequence $\{\omega_n\}\subset K_{\mu}$ satisfying $\omega_n=\omega_{1,n}+\omega_{2,n}$, $\omega_{1,n}$, $\omega_{2,n}\geq 0$, and for $0<\alpha<\mu$,

\begin{align*}
\alpha=\int_{\mathbb{R}^{2}_{+}}x_2\omega_{1,n}\dd x,\quad \mu-\alpha=\int_{\mathbb{R}^{2}_{+}}x_2\omega_{2,n}\dd x, \quad
\textrm{dist}\ (\textrm{spt}\ \omega_{1,n},\textrm{spt}\ \omega_{2,n} )\to \infty.
\end{align*}\\
Observe that for example if $\omega_{1,n}$ and $\omega_{2,n}$ are compactly supported and move away for the $x_1$-direction, the sequence $\{\omega_{n}\}$ is not compact in $L^{2}$. If we have the strict subadditivity of $I_{\mu}$, i.e., $I_{\mu}<I_{\alpha}+I_{\mu-\alpha}$ for $0<\alpha<\mu$, we immediately conclude that this can not occur by sending $n\to\infty$ to $E_{2}[\omega_n]\leq E_{2}[\omega_{1,n}]+E_{2}[\omega_{2,n}]+o(1)$.

The main difficulty is the fact that $K_{\mu}$ has the multiple constraints (impulse $=\mu$, mass $\leq 1$) which is an obstacle to deduce the strict subadditivity of $I_{\mu}$ from the scaling property of $E_{2}$. See \cite[Corollary II.1]{Lions84a}, \cite[p.279]{Lions84b}. We overcome this difficulty by reducing the problem to compactness of a sequence satisfying (1.7) and existence of minimizers of (1.6) by using the Steiner symmetrization $\omega^{*}_{i,n}$, i.e., a rearrangement of $\omega_{i,n}$ satisfying (1.7), $E_{2}[\omega_{i,n}]\leq E_{2}[\omega_{i,n}^{*}]$, conserving $L^{q}$-norms, $1\leq q\leq 2$, and impulse. Since $\omega_{i,n}^{*}$ is non-increasing for $x_1>0$, we are able to show that the weak convergence $\omega_{i,n}^{*}\rightharpoonup \overline{\omega}_i$ in $L^{2}$ implies the convergence of the kinetic energy $E[\omega_{i,n}^{*}]\to E[\overline{\omega}_i]$. This yields

\begin{equation*}
\begin{aligned}
&-I_{\mu}\leq E_{2}[\overline{\omega}_1]+E_{2}[\overline{\omega}_2], \\
&\alpha\geq \int_{\mathbb{R}^{2}_{+}}x_2\overline{\omega}_1\dd x,\quad  \mu-\alpha\geq \int_{\mathbb{R}^{2}_{+}}x_2\overline{\omega}_2\dd x,\quad ||\overline{\omega}_1||_1+||\overline{\omega}_2||_1\leq 1.
\end{aligned}
\end{equation*}\\
A contradiction is deduced from the existence of minimizers of (1.6) (satisfying (1.7)). Indeed, there exists a maximizer $\omega_{1}$ of $E_2$ (a minimizer of $-E_2$) under the constraints $\int x_2\overline{\omega}_1\dd x\leq \alpha$ and $||\overline{\omega}_1||_1\leq 1-||\overline{\omega}_2||_1$ for fixed $\overline{\omega}_2$. The maximizer satisfies $\int x_2\omega_1\dd x=\alpha$ with compact support. Therefore we are able to replace $\overline{\omega}_1$ with $\omega_1$ and apply the same for $\overline{\omega}_2$ for fixed $\omega_1$. Since we can assume that $\textrm{spt}\ \omega_1\cap \textrm{spt}\ \omega_2=\emptyset$ by translation for the $x_1$-variable, 

\begin{align*}
-I_{\mu}\leq E_{2}[\omega_1]+E_{2}[\omega_2]=E_{2}[\omega_1+\omega_2]-\int_{\mathbb{R}^{2}_{+}}\int_{\mathbb{R}^{2}_{+}} G(x,y)\omega_1(x)\omega_2(y)\dd x\dd y\leq -I_{\mu}.
\end{align*}\\
This implies $\omega_{i}\equiv 0$ for $i=1$ or $2$, a contradiction to $\mu=\int x_2(\omega_1+\omega_2)\dd x$.

The existence of the minimizer $\omega_1$ follows from the compactness of a minimizing sequence satisfying (1.7). Since we can assume that a minimizing sequence satisfies (1.7) by the Steiner symmetrization, the existence of the minimizer $\omega_1$ follows from the convergence of the kinetic energy.

\vspace{15pt}
This paper is organized as follows. In Section 2, we prove that $I_{\mu}$ is negative and decreasing for $\mu\in (0,\infty)$ and that minimizers of (1.6) are solutions of (1.4) with compact support. In Section 3, we prove compactness of the kinetic energy for a sequence satisfying (1.7) and existence of minimizers of (1.6). In Section 4, we prove Theorem 1.3 by a concentration compactness principle. In Section 5, we prove existence of symmetric global weak solutions to (1.1) and deduce Theorem 1.4 by a contradiction argument. In Section 5, we prove Theorem 1.5 by the moving plane method.

\vspace{15pt}

\section{A minimization problem}

\vspace{15pt}

We begin with estimates for the kinetic energy $E[\omega]$. Thanks to the finiteness of the impulse $x_2\omega\in L^{1}$, the kinetic energy is finite for $\omega \in L^{2}\cap L^{1}$ and agrees with the Dirichlet energy for the stream function. By using energy estimates, we show that $I_{\mu}$ is decreasing for $\mu\in (0,\infty)$ and any minimizing sequence of $I_{\mu}$ is a bounded sequence in $L^{2}$. In the subsequent section, we prove properties of minimizers.

\vspace{15pt}

\subsection{Properties of $I_{\mu}$}

For the later usage in the proofs of Theorems 1.3 and 1.4, we estimate difference of two energies.

\vspace{15pt}

\begin{prop}
The estimates 

\begin{align*}
&\int_{\mathbb{R}^{2}_{+}}G(x,y)\omega(y)\dd y \leq Cx_2^{1/2}||\omega||^{1/2}_{1}||\omega||^{1/2}_{2}, \tag{2.1} \\
&E[\omega]\leq C||x_2\omega||_{1}^{1/2}{||\omega||_{1}}||\omega||_{2}^{1/2},   \tag{2.2} \\
&\int_{\mathbb{R}^{2}_{+}}\int_{\mathbb{R}^{2}_{+}}G(x,y)\omega_1(x)\omega_2(y)\dd x\dd y\leq C||\omega_1||_{1}^{1/2}||\omega_1||_{2}^{1/2}||x_2\omega_2||_{1}^{1/2}||\omega_2||_1^{1/2},   \tag{2.3}\\
&\left|E[\omega_1]-E[\omega_2]\right|
\leq C||\omega_1-\omega_2||_{1}^{1/2}||\omega_1-\omega_2||_{2}^{1/2}
||x_2(\omega_1+\omega_2)||_{1}^{1/2}||\omega_1+\omega_2||_{1}^{1/2},  \tag{2.4}
\end{align*}\\
hold for $\omega, \omega_i \in L^{2}\cap L^{1}(\mathbb{R}^{2}_{+})$ satisfying $x_2\omega, x_2\omega_i\in L^{1}(\mathbb{R}^{2}_{+})$, $\omega,\omega_i\geq 0$, with some constant $C$, independent of $\omega$, $\omega_{i}$, $i=1,2$.
\end{prop}

\vspace{15pt}

\begin{proof}
The estimate (2.2) follows from (2.3). We suppress the integral region. Observe that  

\begin{align*}
2(E[\omega_1]-E[\omega_2])&=\iint G(x,y)\omega_1(x)\omega_1(y)\dd x\dd y-\iint G(x,y)\omega_2(x)\omega_2(y)\dd x\dd y \\
&=\iint G(x,y)\tilde{\omega}(x)\omega_1(y)\dd x\dd y+\iint G(x,y)\omega_2(x)\tilde{\omega}(y)\dd x\dd y,
\end{align*}\\
for $\tilde{\omega}=\omega_1-\omega_2$ and by $G(x,y)=G(y,x)$, 

\begin{align*}
\iint G(x,y)\omega_2(x)\tilde{\omega}(y)\dd x\dd y=\iint G(y,x)\omega_2(y)\tilde{\omega}(x)\dd x\dd y 
&=\iint G(x,y)\tilde{\omega}(x)\omega_2(y)\dd x\dd y.
\end{align*}\\
We see that

\begin{align*}
2(E[\omega_1]-E[\omega_2])=\iint G(x,y)\tilde{\omega}(x)\hat{\omega}(y)\dd x\dd y,\quad \hat{\omega}=\omega_1+\omega_2.
\end{align*}\\
Thus (2.4) follows from (2.3). 

We set $\psi_1$ by $\omega_1$ and (1.5). By the H\"older's inequality, for $q\in (1,2)$, $1/q=\theta+(1-\theta)/2$, 

\begin{align*}
\psi_1(x)\leq \left(\int_{\mathbb{R}^{2}_{+}}G(x,y)^{q'}\dd y\right)^{1/q'}||\omega_1||_{q}\leq Cx_2^{2/q'}||\omega_1||_{q}
\leq Cx_2^{1-\theta}||\omega_1||_1^{\theta}||\omega_1||_2^{1-\theta}.
\end{align*}\\
Taking $\theta=1/2$ implies (2.1) and 

\begin{align*}
\iint G(x,y)\omega_1(y)\omega_2(x)\dd x\dd y=\int \psi_1(x)\omega_2(x)\dd x&\leq C||\omega_1||_1^{1/2}||\omega_1||_2^{1/2}\int x_2^{1/2}\omega_2(x)\dd x\\
&\leq C||\omega_1||_1^{1/2}||\omega_1||_2^{1/2}||x_2\omega_2||_1^{1/2}||\omega_2||_1^{1/2}.
\end{align*}
We obtained (2.3). This completes the proof.
\end{proof}

\vspace{15pt}

We show that the Dirichlet integral of the stream function is finite.

\vspace{15pt}

\begin{prop}
For $\omega\in L^{2}\cap L^{1}(\mathbb{R}^{2}_{+})$ satisfying $x_2\omega\in L^{1}(\mathbb{R}^{2}_{+})$ and $\omega\geq 0$ ($\omega\nequiv 0$), the stream function (1.5) satisfies $\psi>0$ in $\mathbb{R}^{2}_{+}$, 

\begin{align*}
&\psi(x)\to 0\quad \textrm{as}\ |x|\to\infty,   \tag{2.5} \\
&E[\omega]=\frac{1}{2}||\nabla \psi||_{2}^{2}.   \tag{2.6}
\end{align*}
\end{prop}

\vspace{5pt}

\begin{proof}
By

\begin{align*}
\psi(x)=\int_{\mathbb{R}^{2}_{+}}G(x,y)\omega(y)\dd y=\int_{|x-y|\geq x_2/2}+\int_{|x-y|< x_2/2},
\end{align*}\\
and $G(x,y)\leq \pi^{-1}x_2y_2|x-y|^{-2}$, 

\begin{align*}
\int_{|x-y|\geq x_2/2}G(x,y)\omega(y)\dd y\leq \frac{4}{\pi x_2}||y_2\omega||_1.
\end{align*}\\
By the H\"older's inequality, $1/q+1/q'=1$, $1/q=\theta+(1-\theta)/2$, 

\begin{align*}
\int_{|x-y|< x_2/2}G(x,y)\omega(y)\dd y
&\leq \left(\int_{|x-y|< x_2/2}G(x,y)^{q'}\dd y \right)^{1/q'}\left(\int_{|x-y|< x_2/2}\omega(y)^{q}\dd y \right)^{1/q} \\
&\leq Cx_2^{2/q'}||\omega||_{L^{1}(|x-y|< x_2/2)}^{\theta}||\omega||_{L^{2}(|x-y|< x_2/2)}^{1-\theta}.
\end{align*}\\
Since 

\begin{align*}
\int_{|x-y|< x_2/2}\omega(y)\dd y\leq \frac{2}{x_2}||y_2\omega||_{1},
\end{align*}\\
we have 

\begin{align*}
\int_{|x-y|< x_2/2}G(x,y)\omega(y)\dd y
\leq \frac{C}{x_2^{4/q-3}}||x_2\omega||_{L^{1}}^{\theta}||\omega||_{L^{2}\cap L^{1}}^{1-\theta}.
\end{align*}\\
Hence by (2.1) and for $\delta\in (0,1)$, by taking $q\in (1,2]$ sufficiently small,

\begin{align*}
\psi(x)\leq \frac{C_{\delta}}{(1+x_2)^{1-\delta}}\left(||x_2\omega||_{L^1}+||\omega||_{L^2\cap L^{1}}\right),\quad x\in \mathbb{R}^{2}_{+}.    \tag{2.7}
\end{align*}\\
We take a sequence $\{\omega_n\}\subset C^{\infty}_{c}(\mathbb{R}^{2}_{+})$ such that $\omega_n\to \omega$ in $L^{2}\cap L^{1}(\mathbb{R}^{2}_{+})$ and $x_2\omega_n\to x_2\omega$ in $L^{1}(\mathbb{R}^{2}_{+})$. By (2.7), 

\begin{align*}
\psi(x)
&=\int_{\mathbb{R}^{2}_{+}}G(x,y)(\omega(y)-\omega_n(y))\dd y+\int_{\mathbb{R}^{2}_{+}}G(x,y)\omega_n(y)\dd y \\
&\leq C\left(||x_2(\omega-\omega_n)||_{L^1}+||\omega-\omega_n||_{L^2\cap L^{1}}\right)+\frac{x_2}{\pi \inf_{y\in \textrm{spt}\ \omega_n}|x-y|^{2}}||y_2\omega_n||_{L^{1}}.
\end{align*}\\
Sending $|x|\to\infty$ and then $m\to\infty$ imply (2.5).

We take a non-increasing function $\theta\in C^{\infty}_{c}[0,\infty)$ satisfying $\theta=1$ in $[0,1]$, $\theta=0$ in $[2,\infty)$ and set the cut-off function by $\theta_R(x)=\theta(|x|/R)$. Since $-\Delta \psi=\omega$ in $ \mathbb{R}^{2}_{+}$ and $\psi(x_1,0)=0$, by multiplying $\psi\theta_{R}$ by $-\Delta \psi=\omega$ and integration by parts, 

\begin{align*}
\int_{\mathbb{R}^{2}_{+} }\left(|\nabla \psi|^{2}\theta_R-\frac{1}{2}\psi^{2}\Delta \theta_{R}\right)\dd x=\int_{\mathbb{R}^{2}_{+} } \psi\omega\theta_{R}\dd x. 
\end{align*}\\
Since $\psi\to0$ as $|x|\to\infty$ by (2.5), the second term vanishes as $R\to\infty$. Hence (2.6) follows from the monotone convergence theorem.
\end{proof}

\vspace{15pt}

The function $I_{\mu}$ is negative and decreasing for $\mu\in (0,\infty)$ by (2.2).

\vspace{15pt}

\begin{lem}
\begin{align*}
&I_0=0,  \tag{2.8}\\
&-\infty<I_{\mu}<0,\quad 0<\mu<\infty,  \tag{2.9} \\
&I_{\mu}<I_{\alpha},\quad 0<\alpha <\mu.  \tag{2.10} 
\end{align*}
\end{lem}

\vspace{5pt}

\begin{proof}
Since 

\begin{align*}
I_{\mu}=-\sup_{\omega\in K_{\mu}}E_{2}[\omega],\qquad 
E_2[\omega]=E[\omega]-\frac{1}{2}\int_{\mathbb{R}^{2}_{+}}\omega^{2}\dd x,
\end{align*}\\
we shall show that 

\begin{align*}
&0<\sup_{\omega \in K_{\mu}}E_2[\omega]<\infty,\quad 0<\mu<\infty,  \tag{2.11} \\
&\sup_{\omega \in K_{\alpha}}E_2[\omega]<\sup_{\omega \in K_{\mu}}E_2[\omega],\quad 0<\alpha<\mu.  \tag{2.12}
\end{align*}\\
The property (2.8) is trivial since $K_{0}=\{0\}$. By (2.2) and the Young's inequality, 

\begin{align*}
E_{2}[\omega]\leq C||x_2\omega||_1^{2/3}||\omega||_1^{4/3}\leq C\mu^{2/3},\quad \omega\in K_{\mu}.  
\end{align*} \\
Thus $\sup_{\omega\in K_{\mu}}E_2[\omega]<\infty$. We set $\omega_1=1_{B}$ for $B=B(0,a)$ and choose $a>0$ so that $\int x_2\omega_1\dd x=\mu$. Set $\omega_{\sigma}(x)=\sigma^{3}\omega_1(\sigma x)$, $\sigma>0$, and observe that 

\begin{align*}
&\int_{\mathbb{R}^{2}_{+}}x_2\omega_\sigma\dd x=\int_{\mathbb{R}^{2}_{+}}x_2\omega_1\dd x=\mu, \\
&\int_{\mathbb{R}^{2}_{+}}\omega_\sigma\dd x=\sigma \int_{\mathbb{R}^{2}_{+}}\omega_1\dd x, \\
&E_2[\omega_{\sigma}]=\sigma^{2}\left(E[\omega_1]-\frac{\sigma^{2}}{2}\int_{\mathbb{R}^{2}_{+}}\omega^{2}_1\dd x\right).
\end{align*}\\
Thus for sufficiently small $\sigma>0$, $\omega_{\sigma}\in K_{\mu}$ and 

\begin{align*}
\sup_{\omega\in K_{\mu}}E_{2}[\omega]\geq E_{2}[\omega_{\sigma}]>0.
\end{align*}\\
We proved (2.11). 

It remains to show (2.12). For $\omega\in K_{\alpha}$, $\omega_{\tau}(x)=\tau^{-2}\omega(\tau^{-1}x)$, $\tau>1$, satisfies 

\begin{align*}
&\int_{\mathbb{R}^{2}_{+}}x_2\omega_{\tau}(x)\dd x=\tau\int_{\mathbb{R}^{2}_{+}}x_2\omega(x)\dd x=\tau\alpha, \\
&\int_{\mathbb{R}^{2}_{+}}\omega_{\tau}(x)\dd x=\int_{\mathbb{R}^{2}_{+}}\omega(x)\dd x\leq 1. 
\end{align*}\\
Hence $\omega_{\tau}\in K_{\tau\alpha}$ and 

\begin{align*}
\sup_{\tilde{\omega}\in K_{\tau\alpha}}E_2[\tilde{\omega}]
\geq E_2[\omega_{\tau}]
=E[\omega]-\frac{1}{2\tau^{2}}\int_{\mathbb{R}^{2}_{+}}\omega^{2}\dd x 
=E_{2}[\omega]+\frac{1}{2}\left(1-\frac{1}{\tau^{2}}\right)\int_{\mathbb{R}^{2}_{+}}\omega^{2}\dd x>E_2[\omega]. 
\end{align*}\\
By taking a supremum for $\omega\in K_{\alpha}$, 

\begin{align*}
\sup_{\tilde{\omega}\in K_{\tau\alpha}} E_2[\tilde{\omega}]\geq \sup_{\omega\in K_{\alpha}} E_2[\omega].
\end{align*}\\
If $\sup_{\tilde{\omega}\in K_{\tau\alpha}} E_2[\tilde{\omega}]= \sup_{\omega\in K_{\alpha}} E_2[\omega]$, there exists a maximizing sequence $\{\omega_n\}\subset K_{\alpha}$ such that $E_2[\omega_n]\to \sup_{\omega\in K_\alpha}E_2[\omega]$ and $\omega_n\to0$ in $L^{2}$. By (2.2), $E_{2}[\omega_n]\to 0$. This contradicts (2.11). Hence $\sup_{\tilde{\omega}\in K_{\tau\alpha}} E_2[\tilde{\omega}]>\sup_{\omega\in K_{\alpha}} E_2[\omega]$ and (2.12) holds by taking $\tau=\mu/\alpha$. The proof is complete.
\end{proof}

\vspace{15pt}

\begin{rems}
(i) The strict subadditivity 

\begin{align*}
I_{\mu}<I_{\alpha}+I_{\mu-\alpha},\quad 0<\alpha <\mu, 
\end{align*}\\
is unknown, cf. Lions \cite{Lions84a}.

\noindent
(ii) Any minimizing sequence $\{\omega_n\}$ satisfying $\omega_n\in K_{\mu_n}$, $\mu_n\to \mu$ and $-E_{2}[\omega_n]\to I_{\mu}$ is uniformly bounded in $L^{2}$. Indeed, by (2.2), 

\begin{align*}
||\omega||_2^{2}\leq C\left(||x_2\omega||_1^{2/3}||\omega||_{1}^{4/3}-E_2[\omega ] \right),\qquad \omega\in K_{\mu}.  
\end{align*}\\
By $I_{\mu}<0$, $\limsup_{n\to \infty}||\omega_n||_2\leq C\mu^{1/3}$ follows.
\end{rems}

\vspace{15pt}

\subsection{Properties of minimizers}

\vspace{15pt}

We show that minimizers of (1.6) are solutions to (1.4) for some $W>0$ and $\gamma\geq 0$ with compact support. As noted below in Remarks 2.6 (iii), the flux constant $\gamma$ vanishes if $\mu$ is sufficiently small.

\vspace{15pt}

\begin{prop}
Each minimizer $\omega\in S_{\mu}$ satisfies

\begin{equation*}
\begin{aligned}
&\omega=f (\psi-Wx_2-\gamma),  \\
&\psi(x)=\int_{\mathbb{R}^{2}_{+}}G(x,y)\omega(y)\dd y,
\end{aligned}
\tag{2.13}
\end{equation*}\\
for some constants $W,\gamma\geq 0$, uniquely determined by $\omega$.
\end{prop}

\vspace{15pt}

\begin{proof}
The proof follows from a standard argument, e.g., \cite{FT81}, \cite{Friedman82} for vortex rings. Since $I_{\mu}<0$ by (2.9), minimizers are non-trivial. We take a constant $\delta_0>0$ such that $|\{x\in \mathbb{R}^{2}_{+}\ |\ \omega\geq \delta_0\}|>0$. Here $|E|$ denotes the Lebesgue measure of a set $E\subset \mathbb{R}^{2}_{+}$. We take compactly supported $h_1,h_2\in L^{\infty}(\mathbb{R}^{2}_{+})$ such that $\textrm{spt}\ h_i\subset \{\omega\geq \delta_0\}$, $i=1,2$,  

\begin{align*}
\int_{\mathbb{R}^{2}_{+}}h_1(x)\dd x=1,\quad \int_{\mathbb{R}^{2}_{+}}x_2h_1(x)\dd x=0,\\
\int_{\mathbb{R}^{2}_{+}}h_2(x)\dd x=0,\quad \int_{\mathbb{R}^{2}_{+}}x_2h_2(x)\dd x=1.
\end{align*}\\
We take an arbitrary $\delta\in (0,\delta_0)$ and compactly supported $h\in L^{\infty}(\mathbb{R}^{2}_{+})$ such that $h\geq 0$ on $\{0\leq \omega\leq \delta\}$. We set

\begin{align*}
\eta=h-\left(\int_{\mathbb{R}^{2}_{+}}h \dd x  \right)h_1-\left(\int_{\mathbb{R}^{2}_{+}}x_2 h \dd x \right)h_2
\end{align*}\\
so that $\int \eta\dd x=0$ and $\int x_2\eta\dd x=0$. Observe that $\omega+\varepsilon\eta\geq \delta-\varepsilon||\eta||_{\infty}\geq 0$ on $\{\omega\geq \delta\}$ for small $\varepsilon>0$. Since $\eta=h\geq 0$ on $\{0\leq \omega\leq \delta\}$, $\omega+\varepsilon\eta\geq 0$ on $\{0\leq \omega\leq \delta\}$. Hence $\omega+\varepsilon\eta\in K_{\mu}$. Since $\omega$ is a minimizer of (1.6),

\begin{align*}
0\geq  \frac{\dd}{\dd \varepsilon}E_{2}(\omega+\varepsilon\eta)\Bigg|_{\varepsilon=0}
=\int_{\mathbb{R}^{2}_{+}}\left(\psi-\omega \right)\eta\dd x=:E_{2}'(\omega)\eta. 
\end{align*}\\
By the definition of $\eta$, 

\begin{align*}
E_{2}'(\omega)\eta= E_{2}'(\omega)h-E_{2}'(\omega)h_1\left(\int_{\mathbb{R}^{2}_{+}} h\dd x\right)-E_{2}'(\omega)h_2\left(\int_{\mathbb{R}^{2}_{+}}x_2 h\dd x\right).
\end{align*}\\
By setting $\gamma=E_{2}'(\omega)h_1$ and $W=E_{2}'(\omega)h_2$, 

\begin{align*}
0\geq E_{2}'(\omega)h-\gamma\left(\int_{\mathbb{R}^{2}_{+}} h\dd x\right)-W\left(\int_{\mathbb{R}^{2}_{+}}x_2 h\dd x\right)
=
\int_{\mathbb{R}^{2}_{+}}\left(\psi-Wx_2-\gamma-\omega\right)h\dd x
=\int_{0\leq \omega\leq \delta}+\int_{\omega> \delta}.
\end{align*}\\
We set $\Psi=\psi-Wx_2-\gamma$. Since $h$ is an arbitrary function satisfying $h\geq 0$ on $\{0\leq \omega\leq \delta\}$,

\begin{equation*}
\begin{aligned}
\Psi-\omega&=0\quad \textrm{on}\ \{\omega>\delta\},\\
\Psi-\omega&\leq 0\quad \textrm{on}\ \{0\leq \omega\leq \delta\}. 
\end{aligned}
\tag{2.14}
\end{equation*}\\
Since $\delta>0$ is arbitrary, sending $\delta\to 0$ implies 

\begin{equation*}
\begin{aligned}
\Psi-\omega&=0\quad \textrm{on}\ \{\omega>0\},\\
\Psi&\leq 0\quad \textrm{on}\ \{\omega=0\}.
\end{aligned}
\tag{2.15}
\end{equation*}\\
If $\Psi>0$, $\omega=\Psi$. If $\Psi\leq 0$, $\omega=0$. Thus $\omega=\Psi_+$ and (2.13) holds.

We take a sequence $\{x_n\}$, $x_n={}^{t}(x_{1,n},x_{2,n})$, such that $\omega(x_n)\to 0$ and $x_{n,1}\to\infty$, $x_{n,2}\to 0$. By (2.15), 

\begin{align*}
\limsup_{n\to\infty}\left(\psi(x_n)-Wx_{n,2}-\gamma\right)\leq 0.
\end{align*}\\
Hence $\gamma\geq 0$. By taking an another sequence $\{x_n\}$ such that $\omega(x_n)\to0$ and $x_{n,1}\to0$, $x_{n,2}\to\infty$, $W\geq 0$ follows.

We show uniqueness of $W,\gamma$. Suppose that $\omega$ satisfies (2.13) for $W_{*},\gamma_{*}\geq 0$. Then, $\Psi=\psi-W_{*}x_2-\gamma_{*}$ satisfies (2.14) for $\delta\in (0,\delta_0)$. Hence, 

\begin{align*}
0\geq \int_{\mathbb{R}^{2}_{+}}\left(\Psi-\omega \right)h \dd x
&=\int_{\mathbb{R}^{2}_{+}}\left(\psi-\omega-\gamma_{*}-W_*x_2 \right)h \dd x \\
&=E_{2}'(\omega)h-\gamma_{*}\left(\int_{\mathbb{R}^{2}_{+}} h\dd x\right)-W_{*}\left(\int_{\mathbb{R}^{2}_{+}} x_2h \dd x\right),
\end{align*}\\
for compactly supported $h\in L^{\infty}(\mathbb{R}^{2}_{+})$ satisfying $h\geq 0$ on $\{0\leq \omega\leq \delta\}$. By taking $h=\pm h_1, \pm h_2$, $E'_{2}(\omega)h_1=\gamma_*$,  $E'_{2}(\omega)h_2=W_*$ follow. The proof is complete.
\end{proof}

\vspace{15pt}

\begin{rems}

\noindent
(i) The constant $W$ is positive by the identity \cite[p.1062]{Tu83},

\begin{align*}
W=\left(\frac{1}{2\pi} \int_{\mathbb{R}^{2}_{+}}\int_{\mathbb{R}^{2}_{+}}\frac{x_2+y_2}{|x-y^{*}|^{2}}\omega(x)\omega(y)\dd x\dd y\right) \left(\int_{\mathbb{R}^{2}_{+}}\omega(x)\dd x\right)^{-1},\ y^{*}={}^{t}(y_1,-y_2),   \tag{2.16}
\end{align*}\\
for minimizers $\omega\in S_{\mu}$. The identity (2.16) follows by multiplying $\partial_{x_2}\Psi=\partial_{x_2}\psi-W$ by $\omega$ and integration by parts.

\noindent
(ii) Every minimizer $\omega\in S_{\mu}$ for $\gamma>0$ satisfies 

\begin{align*}
\int_{\mathbb{R}^{2}_{+}}\omega\dd x=1.
\end{align*}\\
Indeed, suppose that $\int\omega\dd x<1$. Then, 

\begin{align*}
\eta=h-\left(\int_{\mathbb{R}^{2}_{+}}x_2h\dd x\right) h_2, 
\end{align*}\\
for $h$ and $h_2$ as in the proof of Proposition 2.5, satisfies 

\begin{align*}
\int_{\mathbb{R}^{2}_{+}} (\omega+\varepsilon \eta)\dd x\leq 1,
\end{align*}\\
for small $\varepsilon>0$ and therefore $\omega+\varepsilon \eta\in K_\mu$. By minimality of $\omega$,

\begin{align*}
\psi-Wx_2-\omega&=0\quad \textrm{on}\ \{\omega>0\},\\
\psi-Wx_2&\leq 0 \quad \textrm{on}\ \{\omega=0\}.
\end{align*}\\
This implies (2.13) for $\gamma=0$, a contradiction to $\gamma>0$.

\noindent 
(iii) If $0<\mu\leq M_1$ for some constant $M_1>0$, every minimizer $\omega \in S_{\mu}$ satisfies

\begin{align*}
\int_{\mathbb{R}^{2}_{+}}\omega \dd x<1.
\end{align*}\\
In particular, $\gamma=0$ by (ii). Indeed, suppose that $\int \omega\dd x=1$. By $\mu=\int_{\mathbb{R}^{2}_{+}}x_2\omega\dd x\geq 2\mu\int_{x_2\geq 2\mu}\omega\dd x$, 

\begin{align*}
\int_{0<x_2<2\mu}\omega\dd x=1-\int_{x_2\geq 2\mu}\omega\dd x\geq \frac{1}{2}.
\end{align*}\\
Observe that by $\omega=\Psi_+\leq \psi$, 

\begin{align*}
\int_{0<x_2<2\mu}\omega\dd x
&\leq \int_{0<x_2<2\mu}\dd x\int_{\mathbb{R}^{2}_{+}}G(x,y)\omega(y)\dd y \\
&=\int_{0<y_2<2\mu}\dd y\int_{\mathbb{R}^{2}_{+}}G(y,x)\omega(x)\dd x\\
&=\int_{\mathbb{R}^{2}_{+}}\omega(x)\dd x\int_{0<y_2<2\mu}G(x,y)\dd y 
=\int_{0<x_2<4\mu}\int_{0<y_2<2\mu}+\int_{x_2\geq 4\mu}\int_{0<y_2<2\mu}.
\end{align*}\\
For $0<x_2<4\mu$, we have 

\begin{align*}
\int_{0<y_2<2\mu}G(x,y)\dd y\leq C\mu^{2}.
\end{align*}\\
In fact, by 

\begin{align*}
\int_{0<y_2< 4\mu}G(x,y)\dd y=\int\limits_{\substack{0< y_2< 4\mu,\\ |x-y|< x_2/2}}+\int\limits_{\substack{0< y_2< 4\mu,\\ |x-y|\geq x_2/2}}.
\end{align*}\\
we estimate 

\begin{align*}
\int\limits_{\substack{0< y_2< 4\mu,\\ |x-y|< x_2/2}}G(x,y)\dd y
\leq \frac{1}{4\pi}\int_{|x-y|< x_2/2}\log\left(1+\frac{4x_2y_2}{|x-y|^{2}}\right)\dd y
&=\frac{x_2^{2}}{4\pi}\int_{|z|< 1/2}\log\left(1+\frac{4(1-z_2)}{|z|^{2}}\right)\dd z \\
&\leq C\mu^{2}.
\end{align*}\\
For $|x-y|\geq x_2/2$, the triangle inequality yields $|x-y^{*}|\leq 5|x-y|$ for $y^{*}={}^{t}(y_1,-y_2)$. By $G(x,y)\leq \pi^{-1}x_2y_2|x-y|^{-2}$, 

\begin{align*}
\int\limits_{\substack{0< y_2< 4\mu,\\ |x-y|\geq  x_2/2}}G(x,y)\dd y
\leq \frac{1}{\pi}\int\limits_{\substack{0< y_2< 4\mu,\\ |x-y|\geq x_2/2}}\frac{x_2y_2}{|x-y|^{2}}\dd y
\leq  \frac{25}{\pi}\int\limits_{\substack{0< y_2< 4\mu,\\ |x-y|\geq x_2/2}}
\frac{x_2y_2}{|x-y^{*}|^{2}}\dd y
\leq C\mu^{2}.
\end{align*}\\
Hence we have the desired estimate. 

For $x_2\geq 4\mu $, by $x_2-y_2\geq x_2/2$, 

\begin{align*}
\int_{0<y_2<2\mu}G(x,y)\dd y\leq \frac{x_2}{\pi}\int_{0<y_2<2\mu}\frac{y_2}{|x-y|^{2}}\dd y
\leq C\mu^{2}.
\end{align*}\\
Hence $1/2\leq \int_{0<x_2<2\mu}\omega\dd x\leq C \mu^{2}\to 0$ as $\mu\to 0$, a contradiction.  
\end{rems}

\vspace{15pt}

The positivity of $W>0$ implies compactness of support for minimizers. We denote by $BUC(\overline{\mathbb{R}^{2}_{+}})$ the space of all bounded uniformly continuous functions in $\overline{\mathbb{R}^{2}_{+}}$ and by $C^{\alpha}(\overline{\mathbb{R}^{2}_{+}})$ the space of all H\"older continuous functions of exponent $0<\alpha<1$ in $\overline{\mathbb{R}^{2}_{+}}$. For an integer $k\geq 0$, $BUC^{k+\alpha}(\overline{\mathbb{R}^{2}_{+}})$ denotes the space of all $\psi\in BUC(\overline{\mathbb{R}^{2}_{+}})$ such that $\partial_{x}^{l}\psi\in BUC(\overline{\mathbb{R}^{2}_{+}})\cap C^{\alpha}(\overline{\mathbb{R}^{2}_{+}})$, for $|l|\leq k$.

\vspace{15pt}

\begin{prop}
For $\omega\in S_{\mu}$, the stream function $(2.13)_2$ satisfies $\psi\in BUC^{2+\alpha}(\overline{\mathbb{R}^{2}_{+}})$, $0<\alpha<1$, $\psi/x_2\in BUC^{1+\alpha}(\overline{\mathbb{R}^{2}_{+}})$ and 

\begin{align*}
\frac{\psi(x)}{x_2}\to 0\quad \textrm{as}\ |x|\to\infty.   \tag{2.17}
\end{align*}
\end{prop}

\vspace{15pt}

\begin{proof}
Since $\omega\in L^{1}\cap L^{2}$, the representation $(2.13)_2$ implies $\nabla^{2}\psi\in L^{q}, q\in (1,2)$ and $\nabla \psi\in L^{p}$, $1/p=1/q-1/2$. By $(2.13)_1$ and (2.5), $\psi$ satisfies 

\begin{equation*}
\begin{aligned}
-\Delta \psi(x)=f(\psi-Wx_2-\gamma)&\quad \textrm{in}\ \mathbb{R}^{2}_{+},\\
\psi=0&\quad \textrm{on}\ \partial\mathbb{R}^{2}_{+},\\
\psi\to 0&\quad \textrm{as}\ |x|\to\infty.
\end{aligned}
\tag{2.18}
\end{equation*}\\
By the Lipschitz continuity of $f$, $\partial_{x}^{l}\psi\in L^{p}_{\textrm{ul}}(\overline{\mathbb{R}^{2}_{+}})$, $|l|= 3$. Here, $L^{p}_{\textrm{ul}}(\overline{\mathbb{R}^{2}_{+}})$ denotes the uniformly local $L^{p}$-space in $\overline{\mathbb{R}^{2}_{+}}$. Hence $\psi\in BUC^{2+\alpha}(\overline{\mathbb{R}^{2}_{+}})$ by the Sobolev embedding. Since $\psi(x_1,0)=0$ and 

\begin{align*}
\frac{\psi(x_1,x_2)}{x_2}=\int_{0}^{1}(\partial_{2}\psi)(x_1,x_2s)\dd s,
\end{align*}\\
$\psi/x_2\in BUC^{1+\alpha}(\overline{\mathbb{R}^{2}_{+}})$ follows. By (2.6) and the Hardy's inequality \cite[2.7.1]{Mazya},

\begin{align*}
\left\|\frac{\psi}{x_2} \right\|_{2}\leq 2\left\|\nabla \psi \right\|_{2}, 
\end{align*}\\
$\psi/x_2\in BUC(\overline{\mathbb{R}^{2}_{+}})\cap L^{2}(\mathbb{R}^{2}_{+})$ and (2.17) follows.
\end{proof}

\vspace{15pt}

\begin{lem}
The support of $\omega \in S_{\mu}$ is compact in $\overline{\mathbb{R}^{2}_{+}}$.
\end{lem}

\vspace{5pt}

\begin{proof}
Since $\textrm{spt}\ \omega=\overline{\{x\in \mathbb{R}^{2}_{+}\ |\ \psi(x)-Wx_2-\gamma > 0  \}}$ for $W>0$ and $\gamma\geq 0$ by $(2.13)_1$ and (2.16), 

\begin{align*}
Wx_2\leq \psi(x),\qquad x\in \textrm{spt}\ \omega.  
\end{align*}\\
Since $\psi/x_2\to 0$ as $|x|\to\infty$ by (2.17), the assertion follows.
\end{proof}

\vspace{15pt}

To prove Theorem 1.5 later in Section 6, we state properties of the associated stream function.

\vspace{15pt}

\begin{lem}
For $\omega\in S_{\mu}$, the stream function $\psi\in BUC^{2+\alpha}(\overline{\mathbb{R}^{2}_{+}})$, $0<\alpha<1$, is a positive solution of (2.18) satisfying $\psi/x_2\in BUC^{1+\alpha}(\overline{\mathbb{R}^{2}_{+}})$, (2.17) and for

\begin{align*}
\Omega=\left\{x\in \mathbb{R}^{2}_{+}\ \middle|\ \psi(x)- Wx_2-\gamma > 0\right\},
\end{align*}\\
$\overline{\Omega}$ is compact in $\overline{\mathbb{R}^{2}_{+}}$. If $0<\mu\leq M_1$, $\gamma=0$, where $M_1$ is the constant as in Remarks 2.6 (iii).
\end{lem}

\vspace{5pt}

\begin{proof}
The assertion follows from Propositions 2.2, 2.7, Lemma 2.8 and Remarks 2.6 (iii).
\end{proof}

\vspace{15pt}

\section{Existence of minimizers}

\vspace{15pt}
We prove existence of minimizers satisfying (1.7) by the Steiner symmetrization. If the minimizing sequence $\{\omega_n\}$ satisfies (1.7), the kinetic energy $E[\omega_n]$ is concentrated on a bounded domain $Q=\{x\in \mathbb{R}^{2}_{+}\ |\  |x_1|< AR,\ x_2< R  \}$ and the weak convergence of the sequence $\{\omega_n\}$ in $L^{2}$ implies the convergence of the energy $E[\omega_n]$. Once we have the convergence of the energy, the existence of minimizers easily follows.

\vspace{15pt}

\begin{prop}[Steiner symmetrization]
For $\omega\geq 0$ satisfying $\omega\in L^{2}\cap L^{1}(\mathbb{R}^{2}_{+})$ and $x_2\omega\in L^{1}(\mathbb{R}^{2}_{+})$, there exists $\omega^{*}\geq 0$ such that 

\begin{equation*}
\begin{aligned}
&\omega^{*}(x_1,x_2)=\omega^{*}(-x_1,x_2), \\
&\omega^{*}(x_1,x_2)\ \textrm{is non-increasing for}\ x_1>0. 
\end{aligned}
\tag{3.1}
\end{equation*}\\
Moreover, 

\begin{align*}
&||\omega^{*}||_{q}=||\omega||_{q}\quad 1\leq q\leq 2,    \\
&||x_2\omega^{*}||_{1}=||x_2\omega||_{1},   \ \\
&E(\omega^{*})\geq E(\omega).  
\end{align*}
\end{prop}

\vspace{5pt}

\begin{proof}
See \cite[Appendix I]{FB74}, \cite[p.1053]{Tu83}.
\end{proof}

\vspace{15pt}

For the later usage in the proof of Theorem 1.3, we state a result for general $0<\mu,\nu<\infty$ with $\lambda=1$. We first find a minimizer of $-E_2$ in a slightly larger space $\tilde{K}_{\mu,\nu}\supset K_{\mu,\nu}$ and then prove that the impulse of this minimizer is exactly $\mu>0$. The goal of this section is to prove:

\vspace{15pt}

\begin{lem}
For $0<\mu,\nu<\infty$, set 

\begin{align*}
\tilde{K}_{\mu,\nu}=\left\{\omega\in L^{2}(\mathbb{R}^{2}_{+})\ \middle|\ \omega\geq 0,\ \int_{\mathbb{R}^{2}_{+}} x_2\omega\dd x\leq \mu,\ \int_{\mathbb{R}^{2}_{+}}\omega\dd x\leq \nu    \right\}.
\end{align*}\\

\noindent
(i) There exists $\omega\in \tilde{K}_{\mu,\nu}$ such that  

\begin{align*}
E_2[\omega]=\sup_{\tilde{\omega}\in \tilde{K}_{\mu,\nu}}E_2[\tilde{\omega}].
\end{align*}

\noindent
(ii) This maximizer $\omega\in \tilde{K}_{\mu,\nu}$ satisfies (1.7), 

\begin{align*}
\int_{\mathbb{R}^{2}_{+}}x_2\omega\dd x=\mu,
\end{align*}\\ 
and is with compact support in $\overline{\mathbb{R}^{2}_{+}}$.
\end{lem}

\vspace{15pt}

The proof of Lemma 3.2 is parallel to the case for vortex rings \cite{FT81}, \cite{Friedman82}. We use the monotonicity $(1.7)_2$ and deduce a decay estimate for the stream function for the $x_1$-variable.

\vspace{15pt}

\begin{prop}
Let $\psi$ be the stream function (1.5) for $\omega\in L^{2}\cap L^{1}(\mathbb{R}^{2}_{+})$ satisfying $x_2\omega\in L^{1}(\mathbb{R}^{2}_{+})$ and $\omega\geq 0$. Assume that (1.7) holds for $\omega$. Then, 

\begin{align*}
\psi(x)\leq C\left( \left(\frac{x_2}{A}\right)^{1/2}||\omega||_1^{1/2}||\omega||_2^{1/2}+\frac{1}{A}||\omega||_1+x_2\left(\frac{A}{x_1}\right)^{2}||x_2\omega||_1\right),\quad x_2\leq \frac{|x_1|}{A}.  \tag{3.2}
\end{align*}\\
The constant $C$ is independent of $\omega$ and $A\geq 1$.
\end{prop}

\vspace{5pt}

\begin{proof}
By replacing $A$ to $A/2$, we prove (3.2) for $x_2\leq 2|x_1|/A$ and $A\geq 2$. We may assume that $x_1>0$. Observe that for a non-increasing function $g(t)\geq 0$ for $t>0$, 

\begin{align*}
\int_{t-t/A}^{t+t/A}g(s)\dd s\leq \frac{4}{A}||g||_{L^{1}(0,\infty)}\quad t>0,\ A\geq 2,
\end{align*}\\
by $tg(t)\leq ||g||_{1}$, $t>0$. Applying this to $\omega$ implies

\begin{align*}
\int_{|x_1-y_1|< x_1/A}\omega(y)\dd y&\leq \frac{4}{A}||\omega||_1.
\end{align*}\\
We set 

\begin{align*}
\psi(x)=\int_{|x-y|< x_2/2}+\int_{|x-y|\geq x_2/2}=:\psi_1+\psi_2.
\end{align*}\\
The conditions $x_2\leq 2x_1/A$ and $|x-y|< x_2/2$ imply $|x_1-y_1|< x_1/A$. By the H\"older's inequality for $1/q=\theta+(1-\theta)/2$, $1/q+1/q'=1$, 

\begin{align*}
\psi_1(x)=\int\limits_{\substack{|x-y|< x_2/2, \\ |x_1-y_1|< x_1/A}}G(x,y)\omega(y)\dd y
&\leq \left(\int_{\mathbb{R}^{2}_{+}}G(x,y)^{q'}\dd y \right)^{1/q'}\left(\int_{|x_1-y_1|< x_1/A}\omega^{q}(y)\dd y \right)^{1/q} \\
&\leq Cx_2^{2/q'}||\omega||_{L^{1}(|x_1-y_1|< x_1/A)}^{\theta}||\omega||_{L^{2}(|x_1-y_1|< x_1/A)}^{1-\theta}.
\end{align*}\\
Taking $\theta=1/2$ yields $\psi_1(x)\leq C(x_2/A)^{1/2}||\omega||_{1}^{1/2}||\omega||_{2}^{1/2}$. We set

\begin{align*}
\psi_{2}(x)=\int\limits_{\substack{|x-y|\geq x_2/2, \\ |x_1-y_1|< x_1/A}}
+\int\limits_{\substack{|x-y|\geq x_2/2, \\ |x_1-y_1| \geq x_1/A}}=: \psi_2^{1}+\psi_{2}^{2}.
\end{align*}\\
By $G(x,y)\leq \pi^{-1}x_2y_2|x-y|^{-2}$, 

\begin{align*}
&\psi_{2}^{1}(x)\leq \frac{1}{\pi}\int\limits_{\substack{|x-y|\geq x_2/2, \\ |x_1-y_1|< x_1/A}} \frac{x_2y_2}{|x-y|^{2}}\omega(y)\dd y
\leq \frac{6}{\pi}\int_{|x_1-y_1|< x_1/A} \omega(y)\dd y
\leq \frac{24}{\pi A}||\omega||_{1}, \\
&\psi_{2}^{2}(x)\leq \frac{1}{\pi}\int\limits_{\substack{|x-y|\geq x_2/2, \\ |x_1-y_1|\geq  x_1/A}} \frac{x_2y_2}{|x-y|^{2}}\omega(y)\dd y 
\leq \frac{x_2}{\pi}\left(\frac{A}{x_1}\right)^{2}||y_2\omega||_{1}.
\end{align*}\\
We obtained (3.2).
\end{proof}

\vspace{15pt}

The stream function estimate (3.2) implies that the kinetic energy $E[\omega]$ is concentrated on a bounded domain $Q=\{x\in \mathbb{R}^{2}_{+}\ |\  |x_1|< AR,\ x_2< R  \}$.

\vspace{15pt}

\begin{prop}
Under the assumption of Proposition 3.3, 

\begin{equation*}
\begin{aligned}
\int_{\mathbb{R}^{2}_{+}\backslash Q}\psi(x)\omega(x)\dd x 
\leq \frac{C}{\min\{A,R\}^{1/2}}\left(  ||\omega||_{L^{1}\cap L^{2}}^{2}+||x_2\omega||_{L^1}^{2}    \right).  
\end{aligned}
\tag{3.3}
\end{equation*}\\
The constant $C$ is independent of $\omega$ and $A, R\geq 1$.
\end{prop}

\vspace{5pt}

\begin{proof}
We decompose

\begin{align*}
\int_{\mathbb{R}^{2}_{+}\backslash Q }\psi(x)\omega(x)\dd x
=\int_{x_2\geq R}+\int_{\substack{x_2<R, \\ |x_1|\geq AR}} ,
\end{align*}\\
and estimate by (2.1) 

\begin{align*}
\int_{x_2\geq R}\psi(x)\omega(x)\dd x\leq C||\omega||_{L^{1}}^{1/2}||\omega||_{L^{2}}^{1/2}\int_{x_2\geq R}x_2^{1/2}\omega\dd x
\leq \frac{C}{R^{1/2}}||\omega||_{L^{1}\cap L^{2}}||x_2\omega||_{L^{1}}.
\end{align*}\\
Since $|x_1|\geq AR$ and $x_2<R$ imply $x_2\leq x_1/A$, applying (3.2) yields 

\begin{align*}
\int\limits_{\substack{x_2<R, \\ |x_1|\geq AR}} \psi(x)\omega(x)\dd x
&\leq C\int\limits_{\substack{x_2<R, \\ |x_1|\geq AR}}\left(  \left(\frac{x_2}{A}  \right)^{1/2}  ||\omega||_{L^{2}\cap L^{1}}+\frac{1}{A}||\omega||_{L^{1}}+x_2\frac{1}{R^{2}}||x_2\omega||_{L^{1}} \right)  \omega(x)\dd x \\
&\leq \frac{C}{\min\{A,R\}^{1/2}}\left( ||\omega||_{L^{1}\cap L^{2}}^{3/2}||x_2\omega||_{L^{1}}^{1/2}+||\omega||_{L^{1}\cap L^{2}}^{2}+||x_2\omega||_{L^{1}}^{2}   \right).
\end{align*}\\
By the Young's inequality, (3.3) follows.
\end{proof}

\vspace{15pt}

Proposition 3.4 implies that the kinetic energy $E[\omega]$ is continuous by the weak continuity in a certain proper subset of $L^{2}$. 

\vspace{15pt}

\begin{lem}
Let $\{\omega_n\}$ be a sequence such that 

\begin{align*}
&\sup_{n\geq 1 }\left\{||\omega_n||_{L^{2}\cap L^{1}}+||x_2\omega_n||_{L^{1}}  \right\}<\infty,\\
&\omega_n\rightharpoonup \omega\quad \textrm{in}\ L^{2}(\mathbb{R}^{2}_{+})\quad \textrm{as}\ n\to\infty.
\end{align*}\\
Assume that each $\omega_n$ satisfies (1.7). Then, 

\begin{align*}
E[\omega_{n}]\to E[\omega]\quad \textrm{as}\ n\to\infty.
\end{align*}
\end{lem}

\vspace{5pt}

\begin{proof}
We decompose the energy into two terms

\begin{align*}
2E[\omega_n]=\int_{\mathbb{R}^{2}_{+}}\psi_n(x)\omega_n(x)\dd x
=\int_{Q}+\int_{\mathbb{R}^{2}_{+}\backslash Q},
\end{align*}\\
and observe that

\begin{align*}
\int_{Q}\psi_n(x)\omega_n(x)\dd x
&=\int_{Q}\omega_{n}(x)\dd x\int_{\mathbb{R}^{2}_{+}}G(x,y)\omega_{n}(y)\dd y 
=\int_{Q}\int_{Q}+\int_{Q}\int_{\mathbb{R}^{2}_{+}\backslash Q}.
\end{align*}\\
By $G(x,y)=G(y,x)$, 

\begin{align*}
\int_{Q}\omega_n(x)\dd x \int_{\mathbb{R}^{2}_{+}\backslash Q}G(x,y)\omega_n(y)\dd y
=\int_{Q}\omega_n(y)\dd y \int_{\mathbb{R}^{2}_{+}\backslash Q}G(x,y)\omega_n(x)\dd x 
\leq \int_{\mathbb{R}^{2}_{+}\backslash Q}\psi_{n}(x)\omega_n(x)\dd x. 
\end{align*}\\
Applying (3.3) yields 

\begin{align*}
\left|2E[\omega_n]-\int_{Q}\int_{Q}G(x,y)\omega_n(x)\omega_n(y)\dd x\dd y\right|
\leq 2\int_{\mathbb{R}^{2}_{+}\backslash Q}\psi_n(x)\omega_n(x)\dd x
\leq \frac{C}{\min\{A,R\}^{1/2}}.
\end{align*}\\
By estimating $E[\omega]$ in the same way, 

\begin{align*}
2\left|E[\omega_n]-E[\omega]\right|
\leq \left|\int_{Q}\int_{Q}G(x,y)\left(\omega(x)\omega(y)-\omega_n(x)\omega_n(y)\right)\dd x\dd y \right|
+\frac{C}{\min\{A,R\}^{1/2}}
\end{align*}\\
Since $G(x,y)\in L^{2}(Q\times Q)$ and $\omega_n(x)\omega_n(y)\rightharpoonup\omega(x)\omega(y)$ in $L^{2}(Q\times Q)$, sending $n\to\infty$ and $A,R\to\infty$ imply the desired result.
\end{proof}

\vspace{15pt}

\begin{proof}[Proof of Lemma 3.2]
By the scaling (1.12), we reduce to the case $0<\mu<\infty$, $\nu=1$ with an abbreviated notation $\tilde{K}_{\mu,1}=\tilde{K}_{\mu}$. Let $\{\omega_n\}\subset \tilde{K}_{\mu}$ be a maximizing sequence of $E_2$. By the Steiner symmetrization, we may assume that $\omega_n$ satisfies (1.7). Since $\{\omega_n\}$ is uniformly bounded in $L^{2}$ as we proved in Remarks 2.4 (ii), by choosing a subsequence (still denoted by $\{\omega_n\}$), there exists $\omega\in L^{2}$ such that $\omega_n\rightharpoonup \omega$ in $L^{2}$ and $||\omega||_2\leq \liminf_{n\to\infty}||\omega_n||_2$. The limit $\omega$ belongs to $\tilde{K}_{\mu}$ and satisfies (1.7). Since $\{\omega_n\}$ satisfies the assumption of Lemma 3.5,

\begin{align*}
\sup_{\tilde{\omega}\in \tilde{K}_{\mu}}E_{2}[\tilde{\omega}]
=\lim_{n\to \infty}E_{2}[\omega_n]
=\lim_{n\to \infty}E[\omega_n]-\frac{1}{2}\liminf_{n\to\infty}||\omega_n||_2^{2}
\leq E[\omega]-\frac{1}{2}||\omega||_2^{2}=E_2[\omega].
\end{align*}\\
Thus $\omega\in \tilde{K}_{\mu}$ is a maximizer. We proved (i).

Since $\sup_{\omega\in \tilde{K}_{\mu}}E_2[\omega]>0$ as we proved (2.9), the maximizer $\omega$ is a non-trivial function and satisfies (2.13) for some constants $W,\gamma\geq 0$ as in Proposition 2.5. By the identity (2.16), we have $W>0$. It remains to show 

\begin{align*}
\int_{\mathbb{R}^{2}_{+}}x_2\omega\dd x=\mu.
\end{align*}\\
Suppose that $\int x_2\omega\dd x<\mu$. Then  

\begin{align*}
\eta=h-\left(\int_{\mathbb{R}^{2}_{+}}h\dd x \right)h_1,
\end{align*}\\
for $h$ and $h_1$ as in the proof of Proposition 2.5, satisfies 

\begin{align*}
\int_{\mathbb{R}^{2}_{+}}x_2(\omega+\varepsilon \eta)\leq \mu,
\end{align*}\\
for small $\varepsilon>0$ and hence $\omega+\varepsilon \eta\in \tilde{K}_{\mu}$. By the maximality of $\omega\in \tilde{K}_{\mu}$, 

\begin{align*}
\psi-\gamma-\omega= 0,\quad \textrm{on}\ \{\omega>0\},\\
\psi-\gamma\leq 0,\quad \textrm{on}\ \{\omega=0\}.
\end{align*}\\
This implies (2.13) for $W=0$, a contradiction to $W>0$ thanks to the uniqueness of $W$ by Proposition 2.5. The compactness of $\textrm{spt}\ \omega$ follows from Lemma 2.8. We proved (ii). 
\end{proof}

\vspace{15pt}

\begin{rem}
It is observed from the proof of Lemma 3.2 that after taking the Steiner symmetrization, $\{\omega_n\}$ satisfies $\lim_{n\to \infty}||\omega_n||_2=||\omega||_2$ and hence $\omega_n\to \omega$ in $L^{2}$. We will see in the next section that any maximizing sequence is relatively compact in $L^{2}$ by translation for the $x_1$-variable without the condition (1.7).
\end{rem}

\vspace{15pt}

\section{Concentrated compactness}

\vspace{15pt}

We prove Theorem 1.3. For a minimizing sequence of (1.6) which does not satisfy the symmetric and non-increasing condition (1.7), Lemma 3.5 can not be directly applied to prove compactness of the sequence. Instead, we apply a concentration compactness principle to get compactness of the minimizing sequence up to translation for the $x_1$-variable. The main difficulty appears when we need to exclude the possibility of dichotomy of the sequence since the strict subadditivity of $I_{\mu}$ is unknown as in Remarks 2.4 (i). To overcome this difficulty, we use the idea from the Steiner symmetrization and reduce the problem to the compactness of a symmetric and non-increasing sequence (Lemma 3.5) and the existence of minimizers of (1.6) (Lemma 3.2).

\vspace{15pt}

\subsection{The case for fixed impulse}

We start with proving Theorem 1.3 for minimizing sequences $\{\omega_n\}\subset K_{\mu}$ of $I_{\mu}$ with fixed impulse. 

\vspace{15pt}

\begin{lem}
Let $0<\mu<\infty$. Let $\{\rho_n\}\subset L^{1}(\mathbb{R}^{2}_{+})$ satisfy 

\begin{align*}
\rho_n\geq 0,\quad \int_{\mathbb{R}^{2}_{+}}\rho_n\dd x=\mu,\quad n\geq 1.
\end{align*}\\
Then, there exists a subsequence $\{\rho_{n_k}\}$ satisfying the one of the followings:

\noindent 
(i) (Compactness)
There exists a sequence $\{y_k\}\subset \overline{\mathbb{R}^{2}_{+}}$ such that $\rho_{n_k}(\cdot +y_k)$ is tight, i.e., for arbitrary $\varepsilon>0$ there exists $R>0$ such that 

\begin{align*}
\int_{B(y_k,R)\cap \mathbb{R}^{2}_{+}}\rho_{n_k}\dd x\geq \mu-\varepsilon,\qquad \textrm{for all}\ k\geq 1.   \tag{4.1}
\end{align*}\\

\noindent 
(ii) (Vanishing) For each $R>0$,

\begin{align*}
\lim_{k\to\infty}\sup_{y\in \mathbb{R}^{2}_{+} }\int_{B(y,R)\cap \mathbb{R}^{2}_{+}}\rho_{n_k}\dd x=0.    \tag{4.2}
\end{align*}\\

\noindent 
(iii) (Dichotomy) There exists $\alpha\in (0,\mu)$ such that for arbitrary $\varepsilon>0$ there exist $k_0\geq 1$ and $\{\rho_{k}^{1}\}$, $\{\rho_{k}^{2}\}\subset L^{1}(\mathbb{R}^{2}_{+})$ such that $\textrm{spt}\ \rho^{1}_{k}\cap \textrm{spt}\ \rho^{2}_{k}=\emptyset$, $0\leq \rho_{k}^{i} \leq \rho_{n_k}$, i=1,2, 

\begin{equation*}
\begin{aligned}
&||\rho_{n_k}-\rho_{k}^{1}-\rho_{k}^{2}||_{L^{1}}+
\left|\int_{\mathbb{R}^{2}_{+}}\rho_{k}^{1}\dd x-\alpha  \right|
+\left|\int_{\mathbb{R}^{2}_{+}}\rho_{k}^{2}\dd x-(\mu-\alpha)  \right|
\leq \varepsilon,\qquad \textrm{for}\ k\geq k_0, \\
&\textrm{dist}\ (\textrm{spt}\ \rho^{1}_{k}, \textrm{spt}\ \rho^{2}_{k})\to \infty\quad \textrm{as}\ k\to\infty.
\end{aligned}
 \tag{4.3}
\end{equation*}
\end{lem}

\vspace{5pt}

\begin{proof}
The assertion is proved in \cite[Lemma I.1]{Lions84a} for the whole space by using the L\'evy's concentration function. The proof also applies to a half space.
\end{proof}

\vspace{15pt}

\begin{rem}
The case (i) is further divided into two cases: (a) $\lim\sup_{k\to\infty} y_{2,k}=\infty$ for $y_k={}^{t}(y_{1,k},y_{2,k})$ and (b) $\sup_{k\geq 1} y_{2,k}<\infty$. In the case (b), we may assume that $y_{2,k}=0$ by replacing $R$. In fact, $B({}^{t}(y_{1,k},0),R' )\supset B(y_k,R)$ for $R'=\sup_{k\geq 1}y_{2,k}+R$. Hence

\begin{align*}
\int_{B({}^{t}(y_{1,k},0),R' )}\rho_{n_k}\dd x\geq \mu-\varepsilon,\quad \textrm{for all}\ k\geq 1.
\end{align*}
\end{rem}

\vspace{15pt}
\begin{lem}
The assertion of Theorem 1.3 holds for minimizing sequences $\{\omega_n\}\subset K_{\mu}$ of $I_{\mu}$ with fixed impulse.  
\end{lem}

\vspace{5pt}

\begin{proof}
Let $\{\omega_n\}\subset K_\mu$ be a minimizing sequence of $I_\mu$. By Remarks 2.4 (ii), $\{\omega_n\}$ is uniformly bounded in $L^{2}$. We set $\rho_n=x_2\omega_n$ and apply Lemma 4.1. Then, for a certain subsequence still denoted by $\{\omega_n\}$, one of the three cases, (iii) Dichotomy, (ii) Vanishing, (i) Compactness, should occur. We shall exclude the first two cases to get compactness of the sequence.\\

\noindent
Case 1.\ \textit{Dichotomy:}\\
There exists some $\alpha\in (0,\mu)$ such that for arbitrary $\varepsilon>0$, there exist $k_0\geq 1$ and $\{\omega_{1,n}\},\{\omega_{2,n}\}\subset L^{1}$ such that $\omega_{3,n}=\omega_{n}-\omega_{1,n}-\omega_{2,n}$ satisfies $\textrm{spt}\ \omega_{1,n}\cap \textrm{spt}\ \omega_{2,n}=\emptyset$, $0\leq \omega_{i,n}\leq \omega_n$, $i=1,2,3$, and 

\begin{align*}
&||x_2\omega_{3,n}||_1+|\alpha_n-\alpha|+|\beta_n-(\mu-\alpha)|\leq \varepsilon,\quad \textrm{for}\ n\geq k_0, \\
&\alpha_n=\int_{\mathbb{R}^{2}_{+}}x_2\omega_{1,n}\dd x,\quad \beta_n=\int_{\mathbb{R}^{2}_{+}}x_2\omega_{2,n}\dd x,\\
&d_n=\textrm{dist}\ (\textrm{spt}\ \omega_{1,n}, \textrm{spt}\ \omega_{2,n})\to\infty\quad \textrm{as}\ n\to\infty.
\end{align*}\\
By choosing a subsequence, we may assume that $\alpha_n\to \overline{\alpha}$ and $\beta_n\to \overline{\beta}$. By suppressing the integral region, we see that 

\begin{align*}
2E[\omega_n]&=\iint G(x,y) \omega_n(x)\omega_n(y)\dd x\dd y \\
&=\iint G(x,y) \omega_{1,n}(x)\omega_{1,n}(y)\dd x\dd y
+\iint G(x,y) \omega_{2,n}(x)\omega_{2,n}(y)\dd x\dd y \\
&+2\iint G(x,y) \omega_{1,n}(x)\omega_{2,n}(y)\dd x\dd y 
+\iint G(x,y) (2\omega_{n}(x)-\omega_{3,n}(x))  \omega_{3,n}(y)\dd x\dd y.
\end{align*}\\
Applying (2.3) implies 

\begin{align*}
\left|\iint G(x,y)(2\omega_{n}(x)-\omega_{3,n}(x))\omega_{3,n}(y)\dd x\dd y\right| 
&\leq C||2\omega_n-\omega_{3,n}||_{1}^{1/2}||2\omega_n-\omega_{3,n}||_{2}^{1/2}||x_2\omega_{3,n}||_{1}^{1/2}||\omega_{3,n}||_{1}^{1/2} \\
&\leq C\varepsilon^{1/2}.
\end{align*}\\
Since $G(x,y)\leq \pi^{-1}x_2y_2|x-y|^{-2}$, 

\begin{align*}
\iint G(x,y)\omega_{1,n}(x)\omega_{2,n}(y)\dd x\dd y
=\iint_{|x-y|\geq d_n} G(x,y)\omega_{1,n}(x)\omega_{2,n}(y)\dd x\dd y
\leq \frac{\mu^{2}}{\pi d_n^{2}}.
\end{align*}\\
Hence 

\begin{align*}
E_2[\omega_n]=E[\omega_n]-\frac{1}{2}\int_{\mathbb{R}^{2}_{+}} \omega_n^{2}\dd x
\leq E_2[\omega_{1,n}]+E_2[\omega_{2,n}]+\frac{\mu^{2}}{\pi d_n^{2}}+C\varepsilon^{1/2}.
\end{align*}\\
We take a Steiner symmetrization $\omega_{i,n}^{*}$ of $\omega_{i,n}$ to see that 

\begin{align*}
&E_2[\omega_n]\leq E_2[\omega_{1,n}^{*}]+E_2[\omega_{2,n}^{*}]+\frac{\mu^{2}}{\pi d_n^{2}}+C\varepsilon^{1/2},\\
&||\omega_{1,n}^{*}||_1+||\omega_{2,n}^{*}||_1\leq 1, \quad 
||\omega_{1,n}^{*}||_2+||\omega_{2,n}^{*}||_2\leq C, \\
&\alpha_n=\int_{\mathbb{R}^{2}_{+}}x_2\omega_{1,n}^{*}\dd x,\quad \beta_n=\int_{\mathbb{R}^{2}_{+}}x_2\omega_{2,n}^{*}\dd x.
\end{align*}\\
By choosing a subsequence (still denoted by $\{\omega_{i,n}^{*}\}$), $\omega_{i,n}^{*} \rightharpoonup \overline{\omega}_i^{\varepsilon}$ in $L^{2}$ and $||\overline{\omega}_i^{\varepsilon}||_2\leq \liminf_{n\to\infty}||\omega_{i,n}^{*}||_2$. Since $\omega_{i,n}^{*}$ is symmetric and non-increasing for $x_1>0$, we apply Lemma 3.5 to get the convergence of the kinetic energy

\begin{align*}
\lim_{n\to\infty}E[\omega_{i,n}^{*}]=E[\overline{\omega}_i^{\varepsilon}],\quad i=1,2.
\end{align*}\\
Sending $n\to\infty$ implies that 

\begin{align*}
&-I_{\mu}\leq E_2[\overline{\omega}_1^{\varepsilon}]+E_2[\overline{\omega}_2^{\varepsilon}]+C\varepsilon^{1/2},\\
&||\overline{\omega}_1^{\varepsilon}||_1+||\overline{\omega}_2^{\varepsilon}||_1\leq 1, \quad 
||\overline{\omega}_1^{\varepsilon}||_2+||\overline{\omega}_2^{\varepsilon}||_2\leq C, \\
&\overline{\alpha}\geq \int_{\mathbb{R}^{2}_{+}}x_2\overline{\omega}_1^{\varepsilon}\dd x,\quad \overline{\beta}\geq \int_{\mathbb{R}^{2}_{+}}x_2\overline{\omega}_2^{\varepsilon}\dd x.
\end{align*}\\
Since $\overline{\omega}_i^{\varepsilon}$ for $\varepsilon>0$ is also symmetric and non-increasing for $x_1>0$, applying the same argument for $\overline{\omega}_i^{\varepsilon}$ and sending $\varepsilon\to0$ implies that $\overline{\omega}_{i}^{\varepsilon}\rightharpoonup \overline{\omega}_i$ in $L^{2}(\mathbb{R}^{2}_{+})$ and 

\begin{align*}
&-I_{\mu}\leq E_2[\overline{\omega}_{1}]+E_2[\overline{\omega}_{2}],\\
&||\overline{\omega}_{1}||_1+||\overline{\omega}_{2}||_1\leq 1, \\
&\alpha\geq \int_{\mathbb{R}^{2}_{+}}x_2\overline{\omega}_{1}\dd x,\quad \mu-\alpha\geq \int_{\mathbb{R}^{2}_{+}}x_2\overline{\omega}_{2}\dd x.
\end{align*}\\
If $\overline{\omega}_1\equiv 0$ and $\overline{\omega}_2\equiv 0$, we have $-I_{\mu}\leq 0$, a contradiction to $I_{\mu}<0$ by (2.9). We may assume that $\overline{\omega}_1\nequiv 0$. We set $\nu_1=1-||\overline{\omega}_2||_1>0$ and apply Lemma 3.2 to take a maximizer $\omega_1\in \tilde{K}_{\alpha,\nu_1}$ of 

\begin{align*}
E_2[\omega_1]=\sup_{\omega\in \tilde{K}_{\alpha,\nu_1}}E_2[\omega],
\end{align*}\\
such that $\int x_2\omega_1\dd x=\alpha$ and $\textrm{spt}\ \omega_1$ is compact in $\overline{\mathbb{R}^{2}_{+}}$. Hence 

\begin{align*}
&-I_{\mu}\leq E_2[\omega_{1}]+E_2[\overline{\omega}_{2}],\\
&||\omega_{1}||_1+||\overline{\omega}_{2}||_1\leq 1, \\
&\alpha= \int_{\mathbb{R}^{2}_{+}}x_2\omega_{1}\dd x,\quad \mu-\alpha\geq \int_{\mathbb{R}^{2}_{+}}x_2\overline{\omega}_{2}\dd x.
\end{align*}\\
If $\overline{\omega}_2\equiv 0$, we have $-I_{\mu}\leq -I_{\alpha}$, a contradiction to $I_{\mu}<I_{\alpha}$ by (2.10). We may assume that $\overline{\omega}_2\nequiv 0$. By setting $\nu_{2}=1-||\omega_1||_1>0$ and taking a maximizer $\omega_2\in \tilde{K}_{\mu-\alpha,\nu_2}$ with compact support in the same way, 

\begin{align*}
&-I_{\mu}\leq E_2[\omega_{1}]+E_2[\omega_{2}],\\
&||\omega_{1}||_1+||\omega_{2}||_1\leq 1, \\
&\alpha= \int_{\mathbb{R}^{2}_{+}}x_2\omega_{1}\dd x,\quad \mu-\alpha= \int_{\mathbb{R}^{2}_{+}}x_2\omega_{2}\dd x.
\end{align*}\\
By translation for the $x_1$-variable, we may assume that $\textrm{spt}\ \omega_1\cap \textrm{spt}\ \omega_2=\emptyset$. Since $\omega_1+\omega_2\in K_{\mu}$, 

\begin{align*}
-I_{\mu}\leq E_2[\omega_1]+E_2[\omega_2]
&=E_2[\omega_1+\omega_2]-\iint G(x,y)\omega_1(x)\omega_2(y)\dd x\dd y \\
&\leq -I_{\mu}-\iint G(x,y)\omega_1(x)\omega_2(y)\dd x\dd y
\leq -I_{\mu}.
\end{align*}\\
Hence, $\omega_i\equiv 0$ for $i=1$ or $2$. This contradicts $\mu=\int_{\mathbb{R}^{2}_{+}}x_2(\omega_1+\omega_2)\dd x$. Thus dichotomy does not occur.\\

\noindent
Case 2.\ \textit{Vanishing:}\\
\begin{align*}
\lim_{n\to\infty}\sup_{y\in \mathbb{R}^{2}_{+}}\int_{B(y,R)\cap \mathbb{R}^{2}_{+}}x_2\omega_n\dd x=0,\quad \textrm{for each}\ R>0.
\end{align*}\\
We shall show that $\lim_{n\to\infty}E[\omega_n]=0$. Since $E_2[\omega_n]\leq E[\omega_n]$, this implies $I_{\mu}\geq 0$, a contradiction to $I_{\mu}<0$.

We set 

\begin{align*}
2E[\omega_n]=\iint G(x,y)\omega_n(x)\omega_n(y)\dd x\dd y 
=\iint_{|x-y|\geq R} +\iint_{|x-y|< R}.
\end{align*}\\
Since $G(x,y)\leq \pi^{-1}x_2y_2|x-y|^{-2}$, 

\begin{align*}
\iint_{|x-y|\geq R} G(x,y)\omega_n(x)\omega_n(y)\dd x\dd y\leq \frac{\mu^{2}}{\pi R^{2}}.
\end{align*}\\
We divide the second term into two terms 

\begin{align*}
\iint_{|x-y|<  R} G(x,y)\omega_n(x)\omega_n(y)\dd x\dd y
=\iint\limits_{\substack{|x-y|< R, \\ G\geq  Rx_2y_2}}
+\iint\limits_{\substack{|x-y|< R, \\ G<  Rx_2y_2}}, 
\end{align*}\\ 
and observe that

\begin{align*}
\iint\limits_{\substack{|x-y|< R, \\ G<  Rx_2y_2}}
G(x,y)\omega_n(x)\omega_n(y)\dd x\dd y
\leq R\mu \left(\sup_{y\in \mathbb{R}^{2}_{+}}\int_{B(y,R)\cap \mathbb{R}^{2}_{+}}x_2\omega_n(x)\dd x\right)\to 0\quad \textrm{as}\ n\to\infty. 
\end{align*}\\
We may assume that $R\geq 1$. The condition $G\geq  Rx_2y_2$ implies $|x-y|\leq R^{-1/2}$. Since $|x-y^{*}|\leq 2x_2+ R^{-1/2}$, $y^{*}={}^{t}(y_1,-y_2)$, 

\begin{align*}
G(x,y)=-\frac{1}{2\pi}\left(\log{|x-y|}-\log{|x-y^{*}|} \right)
\leq \frac{1}{\pi}\left(\left|\log{|x-y|}\right|+x_2 \right), 
\end{align*}\\
\begin{align*}
\left(\int_{|x-y|< R^{-1/2}}G(x,y)^{2}\dd y\right)^{1/2}\leq C(R)(1+x_2),
\end{align*}\\
and $C(R)\to 0$ as $R\to\infty$. Hence

\begin{align*}
\iint\limits_{\substack{|x-y|< R, \\ G\geq Rx_2y_2}}G(x,y)\omega_n(x)\omega_n(y)\dd x\dd y
&\leq \iint_{|x-y|< R^{-1/2}}G(x,y)\omega_n(x)\omega_n(y)\dd x\dd y \\
&\leq ||\omega_n||_2 \int_{\mathbb{R}^{2}_{+}}\omega_n(x)\left(\int_{|x-y|< R^{-1/2}}G(x,y)^{2}\dd y\right)^{1/2}\dd x \\
&\leq C(R)' .
\end{align*}\\
Sending $n\to\infty$, and then $R\to\infty$ implies $\lim_{n\to\infty}E[\omega_n]=0$. Thus vanishing does not occur.\\

\noindent
Case 3. \textit{Compactness:}\\
There exists a sequence $\{y_n\}\subset \overline{\mathbb{R}^{2}_{+}}$ such that for arbitrary $\varepsilon>0$, there exists $R>0$ such that 

\begin{align*}
\int_{B(y_n,R)\cap \mathbb{R}^{2}_{+}}x_2\omega_n\dd x\geq \mu-\varepsilon,\qquad \textrm{for all}\ n\geq 1.
\end{align*}\\
By translation for the $x_1$-variable, we may assume that $y_n={}^{t}(0,y_{2,n})$. Then, there are two cases whether (a) $\limsup_{n\to\infty} y_{2,n}=\infty$ or (b) $\sup_{n\geq 1}y_{2,n}<\infty$. We shall first show that the case (a) does not occur.\\

\noindent 
(a) $\limsup_{n\to\infty} y_{2,n}=\infty$. We may assume that $\lim_{n\to\infty}y_{2,n}=\infty$ by choosing a subsequence. We shall show that $\lim_{n\to\infty}E[\omega_n]=0$. This implies $-I_{\mu}=\lim_{n\to\infty}E_{2}[\omega_n]\leq \lim_{n\to\infty}E[\omega_n]=0$, a contradiction to $I_{\mu}<0$.

We set 
\begin{align*}
2E[\omega_n]=\int_{\mathbb{R}^{2}_{+}} \psi_n\omega_n\dd x=\int_{B(y_n,R)\cap \mathbb{R}^{2}_{+}}+\int_{\mathbb{R}^{2}_{+}\backslash B(y_n,R)}, 
\end{align*}\\
for 
\begin{align*}
\psi_n(x)=\int_{\mathbb{R}^{2}_{+}}G(x,y)\omega_n(y)\dd y.
\end{align*}\\
By (2.1), 

\begin{align*}
\int_{B(y_n,R)\cap \mathbb{R}^{2}_{+}}\psi_n\omega_n\dd x
\leq \left\|\frac{\psi_n}{x_2^{1/2}}\right\|_{\infty}\int_{B(y_n,R)\cap \mathbb{R}^{2}_{+}}x_2^{1/2}\omega_n\dd x
\leq \frac{C\mu}{\left(y_{2,n}-R\right)^{1/2}}\to 0\quad \textrm{as}\ n\to\infty.
\end{align*}\\
By the H\"older's inequality,

\begin{align*}
\int_{\mathbb{R}^{2}_{+}\backslash B(y_n,R)}\psi_n\omega_n\dd x
\leq \left\|\frac{\psi_n}{x_2^{1/2}}\right\|_{\infty}\left(\int_{\mathbb{R}^{2}_{+}\backslash B(y_n,R)}x_2\omega_n\dd x\right)^{1/2}\left(\int_{\mathbb{R}^{2}_{+}\backslash B(y_n,R)}\omega_n\dd x\right)^{1/2} 
\leq C\varepsilon^{1/2}.
\end{align*}\\
Thus sending $n\to\infty$, and then $\varepsilon\to 0$ implies $\lim_{n\to\infty}E[\omega_n]=0$. Thus case (a) does not occur.\\

\noindent 
(b) $\sup_{n\geq }y_{2,n}<\infty$. We may assume that $y_{2,n}=0$ by taking sufficiently large $R>0$ as noted in Remark 4.2, i.e., for $B=B(0,R)$, 

\begin{align*}
\int_{B\cap \mathbb{R}^{2}_{+}}x_2\omega_n\dd x\geq \mu-\varepsilon,\qquad \textrm{for all}\ n\geq 1.
\end{align*}\\
Since $\{\omega_n\}$ is uniformly bounded in $L^{2}$, by choosing a subsequence, $\omega_n\rightharpoonup \omega$ in $L^{2}$ for some $\omega$. By sending $n\to\infty$, 

\begin{align*}
\int_{\mathbb{R}^{2}_{+}}x_2\omega \dd x= \mu.
\end{align*}\\
Hence $\omega\in K_{\mu}$. We shall show that 

\begin{align*}
\lim_{n\to \infty}E[\omega_n]=E[\omega].    \tag{4.4}
\end{align*}\\
This implies that 

\begin{align*}
-I_{\mu}=\lim_{n\to\infty}E_2[\omega_n]\leq \lim_{n\to\infty}E[\omega_n]-\frac{1}{2}\liminf_{n\to\infty}||\omega_n||_{2}^{2}\leq  E_2[\omega]\leq -I_{\mu}.
\end{align*}\\
Hence $\lim_{n\to\infty}||\omega_n||_{2}=||\omega||_2$ and $\omega_n\to \omega$ in $L^{2}$ follows. By 

\begin{align*}
\int_{\mathbb{R}^{2}_{+}}x_2|\omega_n-\omega|\dd x=
\int_{B\cap \mathbb{R}^{2}_{+}}x_2|\omega_n-\omega|\dd x
+\int_{\mathbb{R}^{2}_{+}\backslash B}x_2|\omega_n-\omega|\dd x
\leq C||\omega_n-\omega||_2+2\varepsilon,
\end{align*}\\
sending $n\to\infty$ and then $\varepsilon \to0$ implies $x_2\omega_n\to x_2\omega$ in $L^{1}$. Since $E_2[\omega_n]\to E_2[\omega]$, the limit $\omega\in K_{\mu}$ is a minimizer of $I_{\mu}$.

It remains to show (4.4). We decompose

\begin{align*}
2E[\omega_n]=\int_{\mathbb{R}^{2}_{+}}\psi_n\omega_n\dd x=\int_{B\cap \mathbb{R}^{2}_{+}}+\int_{\mathbb{R}^{2}_{+}\backslash B},
\end{align*}\\
and also 

\begin{align*}
\int_{B\cap \mathbb{R}^{2}_{+}}\psi_n\omega_n\dd x
=\int_{B\cap \mathbb{R}^{2}_{+}}\omega(x)\dd x\int_{\mathbb{R}^{2}_{+}}G(x,y)\omega_n(y)\dd y
=\int_{B\cap \mathbb{R}^{2}_{+}}\int_{B \cap \mathbb{R}^{2}_{+}}
+\int_{B\cap \mathbb{R}^{2}_{+}}\int_{\mathbb{R}^{2}_{+}\backslash B}.
\end{align*}\\
Observe that by $G(x,y)=G(y,x)$, 

\begin{align*}
\int_{B\cap \mathbb{R}^{2}_{+}}\omega_n(x)\dd x\int_{\mathbb{R}^{2}_{+}\backslash B}G(x,y)\omega_n(y)\dd y
&=\int_{B\cap \mathbb{R}^{2}_{+}}\omega_n(y)\dd y\int_{\mathbb{R}^{2}_{+}\backslash B}G(x,y)\omega_n(x)\dd x\\
&\leq \int_{\mathbb{R}^{2}_{+}\backslash B}\omega_n(x)\dd x\int_{ \mathbb{R}^{2}_{+}}G(x,y)\omega_n(y)\dd y\\
&= \int_{\mathbb{R}^{2}_{+}\backslash B}\psi_n(x)\omega_n(x)\dd x. 
\end{align*}\\
Hence

\begin{align*}
\left|2E[\omega_n]-\int_{B\cap \mathbb{R}^{2}_{+}}\int_{B \cap \mathbb{R}^{2}_{+}}G(x,y)\omega_n(x)\omega_n(y)\dd x\dd y\right|
\leq 2\int_{\mathbb{R}^{2}_{+}\backslash B}\psi_n(x)\omega_n(x)\dd x.
\end{align*}\\
By  

\begin{align*}
\int_{\mathbb{R}^{2}_{+}\backslash B}\psi_n(x)\omega_n(x)\dd x\leq \left\|\frac{\psi_n}{x_2^{1/2}}\right\|_{\infty}\left(\int_{\mathbb{R}^{2}_{+}\backslash B}x_2\omega_n\dd x\right)^{1/2}\left(\int_{\mathbb{R}^{2}_{+}\backslash B}\omega_n\dd x\right)^{1/2}\leq C\varepsilon^{1/2}, 
\end{align*}\\
and estimating $E[\omega]$ in the same way, 

\begin{align*}
2\left| E[\omega_n]-E[\omega]\right|
\leq \left|\int_{B\cap \mathbb{R}^{2}_{+}}\int_{B \cap \mathbb{R}^{2}_{+}}G(x,y)\left(\omega_n(x)\omega_n(y)-\omega(x)\omega(y)\right)\dd x\dd y \right|
+C\varepsilon^{1/2}.
\end{align*}\\
Since $G(x,y)\in L^{2}(B\times B)$ and $\omega_n(x)\omega_n(y)\rightharpoonup \omega(x)\omega(y)$ in $L^{2}(B\times B)$, sending $n\to\infty$ and $\varepsilon \to0$ yields $\lim_{n\to\infty}E[\omega_n]=E[\omega]$. The proof is now complete.
\end{proof}

\vspace{15pt} 

\subsection{The case for varying impulse}

We now complete the proof of Theorem 1.3. As used in \cite{Lions84a}, \cite{CL82}, the concentration-compactness lemma (Lemma 4.1) is available even if mass is not exactly the same. See also \cite[Lemma 1]{BNL13}.

\vspace{15pt} 

\begin{lem}
Let $0<\mu<\infty$. For any sequence $\{\rho_n\}\subset L^{1}(\mathbb{R}^{2}_{+})$ satisfying 

\begin{align*}
\rho_n\geq 0\quad n\geq 1,\quad \int_{\mathbb{R}^{2}_{+}}\rho_n\dd x=\mu_n\to \mu\quad \textrm{as}\ n\to\infty,
\end{align*}\\
the assertion of Lemma 4.1 holds by replacing (4.1) with

\begin{align*}
\liminf_{k\to\infty}\int_{B(y_k,R)\cap \mathbb{R}^{2}_{+}}\rho_{n_k}\dd x\geq \mu-\varepsilon,
\end{align*}\\
and $(4.3)_1$ with

\begin{equation*}
\begin{aligned}
\limsup_{k\to\infty}\left\{||\rho_{n_k}-\rho_{k}^{1}-\rho_{k}^{2}||_{L^{1}}+
\left|\int_{\mathbb{R}^{2}_{+}}\rho_{k}^{1}\dd x-\alpha  \right|
+\left|\int_{\mathbb{R}^{2}_{+}}\rho_{k}^{2}\dd x-(\mu-\alpha)  \right|\right\}
\leq \varepsilon.
\end{aligned}
\end{equation*}
\end{lem}

\vspace{5pt} 

\begin{proof}
Applying Lemma 4.1 for $\tilde{\rho}_{n}=\rho_n\mu/\mu_n$ yields the desired result. 
\end{proof}

\vspace{15pt}

\begin{proof}[Proof of Theorem 1.3] 
For a minimizing sequence $\{\omega_n\}$ such that $\omega_n \in K_{\mu_n}$, $\mu_n\to \mu$ and $-E_{2}[\omega_n]\to I_{\mu}$ as $n\to\infty$, we set $\rho_n=x_2\omega_n$ and apply Lemma 4.4. Then the desired result follows the same way as in the proof of Lemma 4.3 without significant modifications.
\end{proof}

\vspace{15pt} 

\section{Orbital stability}

\vspace{15pt} 

We prove Theorem 1.4. We first show existence of global weak solutions of (1.1) satisfying the conservations (1.8). To see this, we recall renormalized solutions of DiPerna-Lions \cite{DL89}. 

\vspace{15pt} 

\subsection{Existence of global weak solutions}

We consider the linear transport equation

\begin{equation*}
\begin{aligned}
\partial_{t}\xi+b\cdot \nabla \xi&=0\quad \textrm{in}\ \mathbb{R}^{2}\times (0,T),\\
\xi(x,0)&=\xi_0\quad \textrm{on}\ \mathbb{R}^{2}\times \{t=0\},
\end{aligned}
\tag{5.1}
\end{equation*}\\
with the divergence-free drift $b$, i.e., $\D\ b=0$, satisfying 

\begin{equation*}
\begin{aligned}
&b\in L^{1}(0,T; W^{1,1}_{\textrm{loc}}(\mathbb{R}^{2}) ),\\
&\frac{b}{1+|x|}\in L^{1}(0,T; L^{1}+ L^{\infty}(\mathbb{R}^{2})  ).
\end{aligned} 
\tag{5.2}
\end{equation*}\\
We denote by $L^{0}$ the set of all measurable functions $f$ such that $|\{|f|>\alpha   \}|<\infty$ for each $\alpha\in (0,\infty)$. We say that $\xi\in L^{\infty}(0,T; L^{0})$ is a renormalized solution of $(5.1)_1$ if $\xi$ satisfies 

\begin{align*}
\partial_t\beta(\xi)+b\cdot \nabla \beta(\xi)=0\quad \textrm{in}\ \mathbb{R}^{2}\times (0,T),  \tag{5.3}
\end{align*}\\
for all $\beta \in C^{1}\cap L^{\infty}(\mathbb{R})$ vanishing near zero, in the sense of distribution. It is proved in \cite[Theorem II. 3]{DL89} under the condition (5.2) that for $\xi_0\in L^{0}$ there exists a unique renormalized solution $\xi\in C([0,T]; L^{0})$ of (5.1) and if $\xi_0\in L^{q}(\mathbb{R}^{2})$, $q\in[1,\infty]$, the renormalized solution satisfies $\xi\in C([0,T]; L^{q}(\mathbb{R}^{2}))$ and 

\begin{align*}
||\xi||_{q}(t)=||\xi_0||_{q}\quad \textrm{for all}\ t\geq 0.  \tag{5.4}
\end{align*}\\
As proved in \cite{LNM06}, every global weak solution of (1.1) for $\zeta_0\in L^{q}\cap L^{1}(\mathbb{R}^{2})$, $q\in (1,\infty)$, is a renormalized solution of (5.1) for $b=k*\zeta$. Thus the conservation $(1.8)_1$ holds for the weak solutions by (5.4).

\vspace{15pt} 

\begin{prop}
For symmetric initial data $\zeta_0\in L^{2}\cap L^{1}(\mathbb{R}^{2})$ such that $x_2\zeta_0 \in L^{1}(\mathbb{R}^{2})$ and $\zeta_0\geq 0$ for $x_2\geq 0$, i.e., $\zeta_0(x_1,x_2)=-\zeta_{0}(x_1,-x_2)$, there exists a symmetric global weak solution $\zeta\in BC([0,\infty); L^{2}\cap L^{1}(\mathbb{R}^{2}))$ of (1.1) such that $x_2\zeta\in BC([0,\infty); L^{1}(\mathbb{R}^{2}) )$, $\zeta\geq 0$ for $x_2\geq 0$, 

\begin{align*}
\int_{0}^{\infty}\int_{\mathbb{R}^{2}}\zeta(\partial_t \varphi+v\cdot \nabla \varphi)\dd x\dd t=-\int_{\mathbb{R}^{2}}\zeta_0(x)\varphi(x,0)\dd x  \tag{5.5}
\end{align*}\\
for $v=k*\zeta$ and all $\varphi\in C^{\infty}_{c}(\mathbb{R}^{2}\times [0,\infty))$. This weak solution $\zeta$ satisfies the conservations (1.8).
\end{prop}

\vspace{5pt} 

\begin{proof}
For smooth and symmetric initial data $\zeta_0\in C^{\infty}_{c}$, there exists a symmetric classical solution $\zeta\in BC([0,\infty); L^{2}\cap L^{1})$ of (1.1) \cite{LNX01}. By the conservations (1.8) and the Biot-Savart law $v=k*\zeta$, the solution satisfies 

\begin{equation*}
\begin{aligned}
&\zeta\in L^{\infty}(0,\infty; L^{2}\cap L^{1}), \\
&x_2 \zeta\in L^{\infty}(0,\infty; L^{1}), \\
&v\in L^{\infty}(0,\infty; L^{p}),\ 2\leq p<\infty, \\
&\nabla v\in L^{\infty}(0,\infty; L^{q}),\ 1<q\leq 2. 
\end{aligned}
\tag{5.6}
\end{equation*}\\
Since $v\cdot \nabla v\in L^{\infty}(0,\infty; L^{r})$, $1< r<2$, by the Euler  equation $\partial_t v+(1+\nabla (-\Delta)^{-1}\D) (v\cdot \nabla v)=0$ and $v=\nabla^{\perp}\phi$, $\nabla^{\perp}={}^{t}(\partial_{x_2},-\partial_{x_1})$, 

\begin{equation*}
\begin{aligned}
&\partial_t v\in L^{\infty}(0,\infty; L^{r}),\ 1< r<2, \\
&\partial_t\phi \in L^{\infty}(0,\infty; L^{s}),\ 2<s<\infty.
\end{aligned}
\tag{5.7}
\end{equation*}\\
The function $v$ satisfies the condition (5.2). Indeed, by $v=k*\zeta$,   $k=k1_{B}+k1_{B^{c}}=k_1+k_2$, $B=B(0,1)$, and the Young's inequality, 

\begin{align*}
||v||_{L^{1}+L^{\infty}}\leq ||k_1*\zeta||_{L^{1}}+||k_2*\zeta||_{L^{\infty}}
\leq (||k_1||_{L^{1}}+||k_2||_{L^{\infty}} )||\zeta||_{L^{1}}.
\end{align*}\\
Hence 

\begin{align*}
v\in L^{\infty}(0,\infty; L^{1}+L^{\infty}).   \tag{5.8}
\end{align*}\\
The existence of a global weak solution of (1.1) satisfying (5.5)-(5.8) for symmetric $\zeta_0\in L^{2}\cap L^{1}$, $x_2\zeta_0\in L^{1}$, $\zeta_0\geq 0$ for $x_2\geq 0$, follows by an approximation of $\zeta_0$ by elements of $C^{\infty}_{c}$, e.g., \cite{MaB}. By the condition (5.8) and the consistency \cite[Theorem II.3 (1)]{DL89}, the constructed global weak solution $\zeta$ is a renormalized solution of (5.1). Hence $\zeta\in BC([0,\infty); L^{2}\cap L^{1} ) $ and $(1.8)_1$ holds.

The conservations $(1.8)_2$ and $(1.8)_3$ follow from the weak form (5.5). To see this, we take a cut-off function $\theta\in C^{\infty}_{c}[0,\infty)$, satisfying $\theta\equiv 1$ in $[0,1]$ and $\theta\equiv 0$ in $[2,\infty)$ and set $\theta_R(x)=\theta(|x|/R)$, $R\geq1$. We set $\varphi=x_2\theta_{R}1_{(0,T)}$ for $T>0$. By approximation of $\varphi$ near $t=T$ and substituting it into (5.5) imply 

\begin{align*}
\int_{0}^{T}\int_{\mathbb{R}^{2}}\zeta v\cdot \nabla (x_2\theta_{R})\dd x\dd t=\int_{\mathbb{R}^{2}}x_2\zeta(x,T)\theta_{R}(x)\dd x-\int_{\mathbb{R}^{2}}x_2\zeta_0(x)\theta_{R}(x)\dd x.
\end{align*}\\
Since 

\begin{align*}
\zeta v\cdot \nabla (x_2\theta_R)=\left(\partial_1\left( \frac{1}{2}\left(|v^{2}|^{2}-|v^{1}|^{2} \right)  \right)-\partial_2( v^{1}v^{2}  ) \right)  \theta_R+\zeta v x_2 \cdot \nabla \theta_R,
\end{align*}\\
sending $R\to\infty$ implies $(1.8)_2$. 

To prove $(1.8)_3$, it suffices to show the conservation of the kinetic energy 

\begin{align*}
\int_{\mathbb{R}^{2}}|v(x,T)|^{2}\dd x=\int_{\mathbb{R}^{2}}|v_0(x)|^{2}\dd x.  \tag{5.9}
\end{align*}\\
Since $2E[\omega]=||v||_2^{2}$ by (2.6), $(1.8)_1$ and (5.9) imply $(1.8)_3$. By $(5.6)$ and $(5.7)_1$, observe that 

\begin{align*}
2\int_{0}^{T}\int_{\mathbb{R}^{2}}v\cdot \partial_t v \dd x\dd t
=\int_{\mathbb{R}^{2}}|v(x,T)|^{2}\dd x-\int_{\mathbb{R}^{2}}|v_0(x)|^{2}\dd x.  \tag{5.10}
\end{align*}\\
By (5.6) and approximation of the test functions in (5.5), we have

\begin{align*}
\int_{0}^{T}\int_{\mathbb{R}^{2}}\zeta(\partial_t \varphi+v\cdot \nabla \varphi)\dd x\dd t=\int_{\mathbb{R}^{2}}\zeta(x,T)\varphi(x,T)\dd x-\int_{\mathbb{R}^{2}}\zeta_0(x)\varphi(x,0)\dd x  
\end{align*}\\
for all $\varphi \in L^{\infty}(\mathbb{R}^{2}\times (0,T))$ satisfying $\nabla \varphi, \partial_t \varphi\in L^{\infty}(0,T; L^{s})$, $2<s<\infty$. By (2.1), (2.5), (5.6) and $(5.7)_2$, substituting $\phi$ into the above and integration by parts yield

\begin{align*}
\int_{0}^{T}\int_{\mathbb{R}^{2}}v\cdot \partial_t v\dd x\dd t=\int_{\mathbb{R}^{2}}|v(x,T)|^{2}\dd x-\int_{\mathbb{R}^{2}}|v_0(x)|^{2}\dd x. 
\end{align*}\\
By (5.10), we obtain (5.9). The proof is complete.
\end{proof}

\vspace{15pt} 

\subsection{An application to stability} 
 
We now apply Theorem 1.3 for: 

\vspace{15pt}

\begin{proof}[Proof of Theorem 1.4]
We give a proof for the case $0<\mu<\infty$, $\nu=\lambda=1$. The proof is also applied to the general case $0<\mu,\nu,\lambda<\infty$ by replacing $K_{\mu}, I_{\mu}, S_{\mu}$ to $K_{\mu,\nu}, I_{\mu,\nu,\lambda}, S_{\mu,\nu,\lambda}$, respectively. Suppose that (1.10) were false. Then there exists $\varepsilon_0>0$ such that for $n\geq 1$, there exist $\zeta_{0,n}\in L^{2}\cap L^{1}$ satisfying $\zeta
_{0,n}\geq 0$, $||\zeta_{0,n}||_1\leq 1$ and $t_n\geq 0$ such that a global weak solution in Proposition 5.1 satisfies 

\begin{align*}
&\inf_{\omega\in S_{\mu}}\left\{||\zeta_{0,n}-\omega||_{2}+||x_2(\zeta_{0,n}-\omega)||_{1}\right\}\leq \frac{1}{n},\\
&\inf_{\omega\in S_{\mu}}\left\{||\zeta_{n}(t_n)-\omega||_{2}
+||x_2(\zeta_{n}(t_n)-\omega)||_{1}\right\}
\geq \varepsilon_0.  
\end{align*}\\
We write $\zeta_{n}=\zeta_n(t_n)$ by suppressing $t_n$. We take ${\omega}_n\in S_{\mu}$ such that $||\zeta_{0,n}-{\omega}_n||_{2}+||x_2(\zeta_{0,n}-\omega)||_{1}\to0$. By (2.4), 

\begin{align*}
\left|E_2[\zeta_{0,n}]+I_{\mu}\right|
=\left|E_2[\zeta_{0,n}]-E_2[{\omega}_{n}]\right|
\to 0\quad \textrm{as}\ n\to\infty.
\end{align*}\\
Thus $\{\zeta_{0,n}\}$ is a minimizing sequence such that $\zeta_{0,n}\in K_{\mu_n}$, $\mu_n=\int x_2\zeta_{0,n}\dd x\to \mu$ and $-E_{2}[\zeta_{0,n}]\to I_{\mu}$ as $n\to\infty$. 

By the conservations (1.8), $\zeta_{n}\in K_{\mu_n}$ and 

\begin{align*}
\left|E_{2}[\zeta_n]+I_{\mu}\right|
=\left|E_2[\zeta_{0,n}]+I_{\mu}\right| \to 0\quad \textrm{as}\ n\to\infty.
\end{align*}\\
Hence $\{\zeta_n\}$ is also a minimizing sequence such that $\zeta_{n}\in K_{\mu_n}$, $\mu_n\to \mu$ and $-E_{2}[\zeta_{n}]\to I_{\mu}$. By Theorem 1.3, there exists a sequence $\{y_{n}\}\subset \partial\mathbb{R}^{2}_{+}$ such that, by choosing a subsequence (still denoted by $\{\zeta_n\}$), there exists $\zeta\in L^{2}\cap L^{1}$ such that 

\begin{align*}
\zeta_n(\cdot +y_n)\to {\zeta} \quad \textrm{in}\ L^{2}(\mathbb{R}^{2}_{+}), \\
x_2 \zeta_n(\cdot +y_n)\to x_2 {\zeta} \quad \textrm{in}\ L^{1}(\mathbb{R}^{2}_{+}), 
\end{align*}\\
and the limit $\zeta\in K_{\mu}$ is a minimizer of $I_{\mu}$, i.e., $\zeta\in S_{\mu}$. Sending $n\to\infty$ implies 

\begin{align*}
0=\inf_{\omega\in S_{\mu}}\left\{||\zeta-\omega||_{2}
+||x_2(\zeta-\omega)||_1\right\}
&=\inf_{\omega\in S_{\mu}}\left(\lim_{n\to\infty}
\left\{||\zeta_n-\omega||_{2}
+||x_2(\zeta_n-\omega)||_1\right\}\right)\\
&\geq \liminf_{n\to\infty}\left( \inf_{\omega\in S_{\mu}}
\left\{||\zeta_n-\omega||_{2}
+||x_2(\zeta_n-\omega)||_1\right\}\right)
\geq \varepsilon_0.
\end{align*}\\
We obtained a contradiction. 
\end{proof}

\vspace{15pt} 

\begin{rems}
(i) It is observed from the above proof that the assertion of Theorem 1.4  holds even if impulse of initial data is merely close to $\mu$, i.e., for $\varepsilon>0$ there exists $\delta>0$ such that for $\zeta_0\in L^{2}\cap L^{1}(\mathbb{R}^{2}_{+})$ satisfying $\zeta_0\geq 0$, $||\zeta_0||_1\leq \nu$ and 

\begin{align*}
\inf_{\omega\in S_{\mu,\nu,\lambda}}||\zeta_0-\omega||_2+\left|\int_{\mathbb{R}^{2}_{+}}x_2\zeta_0\dd x-\mu\right| \leq \delta,  \tag{5.11}
\end{align*}\\
there exists a global weak solution of (1.1) satisfying (1.10).

\noindent
(ii) In \cite{BNL13}, orbital stability by the $L^{2}$-norm is proved if initial data $\zeta_0$ is close to a set of minimizers in the same topology as (5.11). 
\end{rems}

\vspace{15pt} 

 \section{Uniqueness of the Lamb dipole}

\vspace{15pt}
We prove Theorem 1.5. For minimizers $\omega\in S_{\mu}$, the associated stream functions are positive solutions of (2.18) for $W>0$ and $\gamma=0$, provided that $0<\mu\leq M_1$ as in Lemma 2.9. Our goal is to prove that such solutions are only translation of the Lamb dipole (1.3) for $\lambda=1$. 

\vspace{15pt}

\subsection{A decay estimate}

We consider positive solutions $\psi>0$ of the problem: 

\begin{equation*}
\begin{aligned}
-\Delta \psi(x)=f(\psi-Wx_2)&\quad \textrm{in}\ \mathbb{R}^{2}_{+},\\
\psi=0&\quad \textrm{on}\ \partial\mathbb{R}^{2}_{+},\\
\psi\to 0&\quad \textrm{as}\ |x|\to\infty,
\end{aligned}
\tag{6.1}
\end{equation*}\\
for some constant $W>0$. 

\vspace{15pt}

\begin{thm}
Let $\psi\in BUC^{2+\alpha}(\overline{\mathbb{R}^{2}_{+}})$, $0<\alpha<1$, be a positive solution of (6.1) for some $W>0$ such that $\psi/x_2\in BUC^{1+\alpha}(\overline{\mathbb{R}^{2}_{+}})$ and $\psi/x_2\to0$ as $|x|\to\infty$ and for $\Omega=\{x\in \mathbb{R}^{2}_{+}\ |\ \psi(x)- Wx_2> 0\}$, $\overline{\Omega}$ is compact in $\overline{\mathbb{R}^{2}_{+}}$. Then, $\psi(x_1,x_2)=\psi_{L}(x_1+q,x_2)$ for some $q\in \mathbb{R}$, where $\psi_{L}=\Psi_{L}+Wx_2$ and $\Psi_{L}$ is the Lamb dipole (1.3) for $\lambda=1$ and the given $W>0$.
\end{thm}

\vspace{15pt}

We reduce (6.1) to the problem in $\mathbb{R}^{4}$. For $y={}^{t}(y_1,y_2,y_3,y_4)\in \mathbb{R}^{4}$, $y'={}^{t}(y_1,y_2,y_3)$, we set $x_1=y_4$, $x_2=|y'|$ and 

\begin{align*}
\varphi(y)=\frac{\psi(x_1,x_2)}{x_2}.   \tag{6.2}
\end{align*}\\
Since $-\psi$ is non-positive, subharmonic and takes a maximum on $\partial\mathbb{R}^{2}_{+}$, by Hopf's lemma \cite[Chapter 2, Theorem 4]{PW}, $\partial_{x_2}\psi(x_1,0)>0$, $x_1\in \mathbb{R}$. Therefore $\varphi$ is positive in $\mathbb{R}^{4}$. The function $\varphi$ is bounded uniformly continuous and H\"older continuous up to first orders in $\mathbb{R}^{4}$, i.e. $\varphi\in BUC^{1+\alpha}(\mathbb{R}^{4})$, $0<\alpha<1$. Moreover, $\varphi$ is continuously differentiable up to second orders in $\mathbb{R}^{4}\backslash \{y'=0\}$ and satisfies 

\begin{align*}
-\Delta_{y}\varphi=-\left(\partial_{x_2}^{2}+\frac{2}{x_2}\partial_{x_2}+\partial_{x_1}^{2}\right)\frac{\psi(x_1,x_2)}{x_2}=f(\varphi-W)\quad \textrm{in}\ \mathbb{R}^{4}\backslash \{y'=0\}.
\end{align*}\\
The function $\varphi$ is regular up to $y'=0$ and satisfies the equation in $\mathbb{R}^{4}$. In fact, by the boundedness of $\nabla \varphi$ near $y'=0$, $\varphi$ satisfies the Poisson equation in $\mathbb{R}^{4}$ in a weak sense. Hence, by $\partial_{y}^{l}f(\varphi-W)\in L^{p}_{\textrm{ul}}(\mathbb{R}^{4})$, $|l|\leq 1$, $1<p<\infty$, and a regularity result for weak solutions \cite{GT}, $\partial_{y}^{l} \varphi\in L^{p}_{\textrm{ul}}(\mathbb{R}^{4})$, $|l|=3$, follows. In particular, $\varphi\in BUC^{2+\alpha}(\mathbb{R}^{4})$, $0<\alpha<1$, by the Sobolev embedding. Hence $\varphi \in BUC^{2+\alpha}(\mathbb{R}^{4})$ is a positive solution of  

\begin{equation*}
\begin{aligned}
-\Delta_y \varphi=f(\varphi-W)&\quad \textrm{in}\ \mathbb{R}^{4},\\
\varphi\to 0&\quad \textrm{as}\ |y|\to\infty.
\end{aligned}
\tag{6.3}
\end{equation*}\\
We set the support of $f(\varphi-W)$ by $\overline{\Xi}$ for 

\begin{align*}
\Xi=\left\{y\in \mathbb{R}^{4}\ \middle|\ \varphi(y)-W> 0\  \right\}.  \tag{6.4}
\end{align*}\\
Since $\overline{\Omega}$ is compact in $\overline{\mathbb{R}^{2}_{+}}$ and $\varphi(y',y_4)=\psi(y_4, |y'|)/|y'|$, $\overline{\Xi}$ is compact in $\mathbb{R}^{4}$ and axisymmetric for the axis $y'=0$.  

\vspace{15pt}

Since $\varphi$ is a positive solution of (6.3), applying the result of \cite[Theorem 4, 2.3. Remark 1]{GNN} implies that $\varphi$ is radially symmetric for some point in $\mathbb{R}^{4}$. See also \cite[Theorem 3.3]{Fra96}, \cite[Theorem 2.1]{Burton96}. We give the proof below for completeness.

\vspace{15pt}

\begin{lem}
There exists $p>0$ and $q\in \mathbb{R}$ such that 

\begin{equation*}
\begin{aligned}
&\varphi(y',y_4+q)=\frac{p}{|y|^{2}}+g(y),\\
&|g(y)|\leq \frac{C}{|y|^{4}},\quad |\nabla g(y)|\leq \frac{C}{|y|^{5}}, \quad \textrm{for}\ |y|\geq 2R+|q|,
\end{aligned}
\tag{6.5}
\end{equation*}\\
for some $R>0$ such that $\Xi\subset B(0,R)$ with some constant $C$, where $B(0,R)$ is an open ball in $\mathbb{R}^{4}$.
\end{lem}

\vspace{5pt}

\begin{proof}
We represent $\varphi$ by the Newton potential of $f(\Phi)$, $\Phi=\varphi-W$, by using the fundamental solution of the Laplace equation in $\mathbb{R}^{4}$, i.e., $\Gamma(y)=(4\pi^{2})^{-1}|y|^{-2}$. By compactness of the support $\Xi\subset \mathbb{R}^{4}$ and $\varphi(y)\to0$ as $|y|\to\infty$, we have 

\begin{align*}
\varphi(y)=\int_{\Xi}\Gamma(y-z)f(\Phi)\dd z.
\end{align*}\\
This implies the expansion

\begin{align*}
&\varphi(y)=\Gamma(y)\int_{\Xi}f(\Phi)\dd z-\nabla_y\Gamma(y)\cdot \left(\int_{\Xi}zf(\Phi)\dd z\right) 
+g_0(y), \\
&|g_0(y)|\leq \frac{C}{|y|^{4}},\quad |\nabla g_0(y)|\leq \frac{C}{|y|^{5}},\quad \textrm{for}\ |y|\geq 2R.
\end{align*}\\
Hence 

\begin{align*}
&\varphi(y)=\frac{p}{|y|^{2}}+\sum_{j=1}^{4}\frac{p_j y_j}{|y|^{4}}+g_0(y), \\
&p=\frac{1}{4\pi^{2}}\int_{\Xi}f(\Phi)\dd z,\quad p_j=\frac{1}{2\pi^{2}}\int_{\Xi}z_jf(\Phi)\dd z,\quad j=1,2,3,4.
\end{align*}\\
Since $\Xi$ and $\Phi$ are symmetric for $y'=0$, $p_j=0$ for $j=1,2,3$. By taking $q=p_4/(2p)$, (6.5) follows.
\end{proof}

\vspace{15pt}

\subsection{Moving plane method}
We apply the moving plane method. The following Propositions 6.3-6.6 are due to \cite[Lemmas 4.1-4.4]{GNN} (see also Lemmas 3.5-3.8 and C.1 of \cite{AF86}). Propositions 6.3, 6.4 are based on the decay estimate (6.5). Proposition 6.5 is by a maximum principle for the monotone function $f$. \\

We take an arbitrary unit vector $n$ in $\mathbb{R}^{4}$ and consider the hyperplane $T_{\kappa}=\{y\in \mathbb{R}^{4}\ |\ y\cdot n=\kappa\ \}$ for $\kappa>0$. By rotation of (6.3), we shall suppose that $n={}^{t}(1,0,0,0)$ and $T_{\kappa}=\{y_1=\kappa\}$. For $y=(y_1,\tilde{y})$, $\tilde{y}=(y_2,y_3,y_4)$, we denote by $y^{\kappa}=(2\kappa-y_1,\tilde{y})$ the reflection with respect to the hyperplane $T_{\kappa}$.

\vspace{15pt}

\begin{prop}
Let $\phi(y)=\varphi(y',y_4+q)$ as in (6.5). Let $\kappa>0$. Consider two points  $y=(y_1,\tilde{y})$ and $z=(z_1,\tilde{y})$ in $\mathbb{R}^{4}$ such that $y_1<z_1$ and $(y_1+z_1)/2\geq \kappa$. There exists $R_{\kappa}>0$ depending only on $\min\{1,\kappa\}$ such that  

\begin{align*}
\phi(y)>\phi(z),\quad \textrm{for}\quad |y|\geq R_{\kappa}.  
\end{align*}
\end{prop}

\vspace{5pt}

\begin{prop}
There exists $\kappa_0\geq 1$ such that for $\kappa\geq \kappa_0$, 

\begin{align*}
\phi(y)>\phi(y^{\kappa}),\quad \textrm{for}\ y_1<\kappa.  
\end{align*}
\end{prop}

\vspace{5pt}

\begin{prop}
Suppose that there exists $\kappa>0$ such that 

\begin{align*}
&\phi(y)\geq \phi(y^{\kappa}),\quad \textrm{for}\ y_1<\kappa,\\
&\phi(y_0)\neq \phi(y_0^{\kappa}),\quad \textrm{for some}\ y_0\in \mathbb{R}^{4}.
\end{align*}\\
Then, $\phi(y)>\phi(y^{\kappa})$ for $y_1<\kappa$ and $\partial_{y_1}\phi(\kappa,\tilde{y})<0$ for $\tilde{y}\in \mathbb{R}^{3}$. 
\end{prop}

\vspace{5pt}

\begin{prop}
The set $\{\kappa>0\ |\ \phi(y)>\phi(y^{\kappa})\ \textrm{for}\ y_1<\kappa \}$ is open in $\mathbb{R}$.
\end{prop}

\vspace{5pt}

\begin{lem}
The function $\phi(y)=\varphi(y', y_4+q)$ is radially symmetric in $\mathbb{R}^{4}$ and decreasing in radial direction.
\end{lem}

\vspace{5pt}

\begin{proof}
By Proposition 6.4, $\phi(y)>\phi(y^{\kappa})$, $y_1<\kappa$, for all $\kappa\geq \kappa_0$ for some $\kappa_0\geq 1$. Since $\{\kappa>0\ |\ \phi(y)>\phi(y^{\kappa})\ \textrm{for}\ y_1<\kappa \}$ is open by Proposition 6.6, we take a maximal interval $(\kappa_*,\infty )$ such that $\phi(y)>\phi(y^{\kappa})$, $y_1<\kappa$, holds for $\kappa>\kappa_{*}$. By continuity, we have $\phi(y)\geq \phi(y^{\kappa_{*}})$ for $y_1<\kappa_{*}$. We shall show that $\kappa_{*}=0$. Suppose that $\kappa_{*}>0$. By Proposition 6.3, there exits $y_0\in \mathbb{R}^{4}$ such that $\phi(y_0)\neq \phi(y_0^{\kappa_{*}})$. By Proposition 6.5, $\phi(y)>\phi(y^{\kappa_*})$ for $y_1<\kappa_{*}$. This contradicts the maximality of $\kappa_{*}$. We thus conclude that $\kappa_{*}=0$ and $\phi(y_1,\tilde{y})\geq  \phi(-y_1,\tilde{y})$ for $y_1\leq 0$. Since $\partial_{y_1}\phi(y_1,\tilde{y})<0$ by Proposition 6.5, $\phi$ is decreasing for $y_1>0$.

Applying the same argument for $n={}^{t}(-1,0,0,0)$ implies that $\phi$ is an even function for $y_1$. By rotation, $\phi$ is symmetric for every unit vectors in $\mathbb{R}^{4}$. Hence $\phi$ is radially symmetric and decreasing.
\end{proof}

\vspace{15pt}

\begin{proof}[Proof of Theorem 6.1]
Since $\phi(y)=\varphi(y', y_4+q)$ is radially symmetric and $|y|=|x|$, $\phi(y)=\phi(|y|)$ and 

\begin{align*}
\frac{\psi(x_1+q,x_2)}{x_2}=\varphi(y',y_4+q)
=\phi(y',y_4)
=\phi\left(|x|\right).
\end{align*}\\
By translation of $\psi$ for the $x_1$-variable, we may assume that $q=0$, i.e., $\psi(x_1,x_2)/x_2=\phi(|x|)$. By the polar coordinate $x_1=r\cos \theta$, $x_2=r\sin\theta$, we set 

\begin{align*}
\Psi(x)=\psi(x)-Wx_2=(\phi(r)-W)r \sin\theta=:\eta(r)\sin\theta.
\end{align*}\\
We prove $\Psi=\Psi_L$. By (6.1), $\Psi$ satisfies 

\begin{equation*}
\begin{aligned}
-\Delta \Psi&=\Psi\quad \textrm{in}\ \Omega,\\
-\Delta \Psi&=0\quad \textrm{in}\ \mathbb{R}^{2}_{+}\backslash\Omega,\\
\Psi&=0\quad \textrm{on}\ \partial\mathbb{R}^{2}_{+}\cup\partial\Omega, \\
\partial_{x_1}\Psi&\to 0,\ \partial_{x_2} \Psi\to -W\quad \textrm{as}\ |x|\to\infty.
\end{aligned}
\tag{6.6}
\end{equation*}\\
Since $\phi(r)$ is decreasing for $r>0$ and $\Psi=0$ on $\partial\Omega$, there exists some $a>0$ such that $\phi(a)=W$ and $\Omega=B(0,a)\cap \mathbb{R}^{2}_{+}$. Substituting $\Psi=\eta(r)\sin\theta$ into $(6.6)_1$ implies that $\eta(r)$ is a solution of the Bessel's differential equation:

\begin{equation*}
\begin{aligned}
\ddot{\eta}+\frac{1}{r}\dot{\eta}-\frac{1}{r^{2}}\eta+\eta=0,\ \eta&>0,\quad 0<r<a,\\
\eta(a)&=0.
\end{aligned}
\tag{6.7}
\end{equation*}\\
Solutions of (6.7) are given by a linear combination of the Bessel functions of the first and second kind of order one. Since $\eta(r)>0$ is bounded at $r=0$ and $\eta(a)=0$, 

\begin{align*}
&\eta(r)=C_1J_1(r),\\
&a=c_0, 
\end{align*}\\
for some constant $C_1$, where $c_0$ is the first zero point of $J_1$. Hence, $\Psi(x)=C_1J_1(r)\sin\theta$ for $r\leq a$.

In a similar way, we consider the region $r\geq a$. Since $\Psi$ is harmonic for $r>a$, $\eta=C_2/r+C_3r$ with some constants $C_2$, $C_3$. Since $\nabla \Psi=(C_2/r^{2}){}^{t}(-\sin{2\theta}, \cos{2\theta})+
{}^{t}(0,C_3)$, sending $r\to\infty$ implies that $C_3=-W$. By $\Psi=0$ for $r=a$, $C_2=Wa^{2}$. Hence $\Psi(x)=-W(r-a^{2}/r)\sin\theta$ for $r>a$.

The constant $C_1$ is determined by continuity of $\partial_r\Psi$ at $r=a$, i.e., $\lim_{r\to a+0}\partial_r\Psi=\lim_{r\to a-0}\partial_r\Psi$. By using $\dot{J}_1(c_0)=J_0(c_0)$, $C_1=-2W/J_0(c_0)=C_L$ follows. We proved $\Psi=\Psi_L$.  
\end{proof}

\vspace{15pt}

\begin{proof}[Proof of Theorem 1.5]
By the scaling (1.12), we reduce to the case $\nu=\lambda=1$. By Theorem 1.3, $S_{\mu}$ is not empty, i.e., $S_{\mu} \neq \emptyset$. Let $0<\mu\leq M_1$ for the constant $M_1>0$ as in Remarks 2.6 (iii). For an arbitrary $\omega\in S_{\mu}$, the associated stream function $\psi$ is a positive solution of (6.1) for some $W>0$ satisfying $\psi/x_2\to0$ as $|x|\to\infty$ and for $\Omega=\{\psi- Wx_2> 0 \}$, $\overline{\Omega}$ is compact in $\overline{\mathbb{R}^{2}_{+}}$ by Lemma 2.9. Applying Theorem 6.1 and $\omega\in K_{\mu}$ imply that $ \omega$ is translation of the Lamb dipole $\omega_L$ for $W=\mu/(c_0^{2} \pi)$. Hence $S_{\mu}\subset \{\omega_{L}(\cdot+y)\ |\ y\in \partial\mathbb{R}^{2}_{+} \}$.

Since $S_{\mu}\neq \emptyset$, there exists $\omega\in S_{\mu}$ and $y_0\in \partial\mathbb{R}^{2}_{+}$ such that $\omega=\omega_{L}(\cdot +y_0)$ for the Lamb dipole $\omega_{L}$ for $W=\mu/(c_0^{2} \pi)$. By translation invariance of $E_2$ for the $x_1$-variable, $\{\omega_{L}(\cdot+y)\ |\ y\in \partial\mathbb{R}^{2}_{+} \}\subset S_{\mu}$ follows. We proved (1.11). The proof is now complete.
\end{proof}



\vspace{15pt}

\section*{Acknowledgements}
The work of the first author is partially supported by JSPS through the Grant-in-aid for Young Scientist (B) 17K14217, Scientific Research (B) 17H02853 and Osaka City University Advanced Mathematical Institute (MEXT Joint Usage / Research Center on Mathematics and Theoretical Physics). The work of the second author is partially supported by NRF-2018R1D1A1B07043065, the Research Fund (1.190136.01) of UNIST (Ulsan National Institute of Science \& Technology) and by the POSCO Science Fellowship of POSCO TJ Park Foundation.

\vspace{15pt}

\bibliographystyle{abbrv}
\bibliography{ref}

\end{document}